
\input amssym

\magnification=\magstep1
\hsize=15,4truecm
\vsize=22truecm
\advance\voffset by 1truecm
\mathsurround=1pt

\def\chapter#1{\par\bigbreak \centerline{\bf #1}\medskip}

\def\section#1{\par\bigbreak {\bf #1}\nobreak\enspace}

\def\sqr#1#2{{\vcenter{\hrule height.#2pt
      \hbox{\vrule width.#2pt height#1pt \kern#1pt
         \vrule width.#2pt}
       \hrule height.#2pt}}}

\def\k{\kappa}
\def\o{\omega}

\def\d{\delta}

\def\l{\lambda}

\def\a{\alpha}
\def\b{\beta}

\def\g{\gamma}

\def\x{\xi}

\def\n{\eta}

\def\tp{\hbox{\rm t}}


\def\A{{\cal A}}
\def\B{{\cal B}}
\def\C{{\cal C}}
\def\D{{\cal D}}

\def\M{{\bf M}}


\def\th #1 #2. #3\par\par{\medbreak{\bf#1 #2.
\enspace}{\sl#3\par}\par\medbreak}
\def\rem #1 #2. #3\par{\medbreak{\bf #1 #2.
\enspace}{#3}\par\medbreak}
\def\proof{{\bf Proof}.\enspace}
\def\sqr#1#2{{\vcenter{\hrule height.#2pt
      \hbox{\vrule width.#2pt height#1pt \kern#1pt
         \vrule width.#2pt}
       \hrule height.#2pt}}}
\def\eop{\mathchoice\sqr34\sqr34\sqr{2.1}3\sqr{1.5}3}

                                                                     %
                                                                     %
\newdimen\refindent\newdimen\plusindent                              %
\newdimen\refskip\newdimen\tempindent                                %
\newdimen\extraindent                                                %
                                                                     %
                                                                     %
\def\ref#1 #2\par{\setbox0=\hbox{#1}\refindent=\wd0                  %
\plusindent=\refskip                                                 %
\extraindent=\refskip                                                %
\advance\extraindent by 30pt                                         %
\advance\plusindent by -\refindent\tempindent=\parindent 
\parindent=0pt\par\hangindent\extraindent 
{#1\hskip\plusindent #2}\parindent=\tempindent}                      %
\refskip=\parindent                                                  %
                                                                     %

\def\empty{\emptyset}

\def\raj{\restriction}

\def\da{\downarrow}

\def\nda{\mathrel{\lower0pt\hbox to 3pt{\kern3pt$\not$\hss}\downarrow}}
\def\nDa{\mathrel{\lower0pt\hbox to 3pt{\kern3pt$\not$\hss}\Downarrow}}
\def\nbot{\mathrel{\lower0pt\hbox to 4pt{\kern3pt$\not$\hss}\bot}}
\def\ekom{\mathrel{\lower3pt\hbox to 0pt{\kern3pt$\sim$\hss}\mapsto}}
\def\do{\triangleright}

\def\anR{\mathrel{\lower1pt\hbox to 2pt{\kern3pt$R$\hss}\not}}
\def\anoR{\mathrel{\lower1pt\hbox to 2pt{\kern3pt$\overline{R}$\hss}\not}}

\def\anRm{\mathrel{\lower1pt\hbox to 2pt{\kern3pt$R^{-}$\hss}\not}}

\def\ndda{\mathrel{\lower0pt\hbox to 1pt{\kern3pt$\not$\hss}\downdownarrows}}

\null
\vskip 2truecm
\centerline{\bf STRONG SPLITTING IN STABLE HOMOGENEOUS MODELS}
\vskip 1truecm
\centerline{Tapani Hyttinen$^{*}$ and Saharon Shelah$^{\dagger}$}
\vskip 2.5truecm

\chapter{Abstract}
\bigskip

In this paper we study elementary submodels
of a stable homogeneous structure.
We improve the independence relation defined in
[Hy]. We apply this to prove a structure theorem.
We also show that dop and sdop
are essentially equivalent, where the negation of dop is the property
we use in our structure theorem and sdop implies nonstructure, see [Hy].

\vskip 2.5truecm

\chapter{1. Basic definitions and spectrum of stability}

The purpose of this paper is to develop theory of independence
for elementary submodels of a homogeneous structure.
We get a model class of this kind if in addition to
it's first-order theory we require that the models omit
some (reasonable) set of types, see [Sh1].
If the set is empty, then we are in the 'classical situation'
from [Sh2]. In other words,
we study
stability theory without the compactness theorem.
So e.g. the theory of $\Delta$-ranks is lost and so
we do not get an independence notion from ranks.
Our independence notion is based on strong splitting.
It satisfies
the basic properties of forking in a rather weak form.
The main problem is finding free extensions. So the arguments
are often based on the definition of the independence notion instead
of the 'independence-calculus'.

\vskip 1.7truecm

\noindent
$*$ Partially supported by the Academi of Finland.

\noindent
$\dagger$ Research supported by the United States-Israel Binational
Science Foundation. Publ. 629.

\vfill
\eject

Throughout this paper we assume that $\M$ is a homogeneous model
of similarity type (=language) $L$
and that $\M$ is $\xi$-stable for some $\xi <\vert\M\vert$
(see [Sh2] I Definition 2.2).
Let $\l (\M )$ be the least such $\xi$. By [Sh1],
$\l (\M )<\beth ((2^{\vert L\vert +\o})^{+})$.
We use $\M$ as a monster model and so we assume that the
cardinality of $\M$ is large enough
for all constructions we
do in this paper. In fact we assume that $\vert\M\vert$
is strongly inaccessible. Alternatively we could
assume less about $\vert\M\vert$ and instead of studying all
elementary submodels of $\M$, we could study suitably small
ones.

Notice $Th(\M )$ may well be unstable.
Notice also
that if $\Delta$ is a stable finite diagram, then $\Delta$
has a monster model like $\M$, see [Sh1].

By a model we mean an elementary submodel of $\M$ of cardinality
$<\vert\M\vert$, we write $\A$, $\B$ and so on for these. So
if $\A\subseteq\B$ are models, then $\A$ is an elementary submodel of
$\B$
Similarly by a set we mean a subset of $\M$ of cardinality
$<\vert\M\vert$,
unless we explicitly say otherwise.
We write $A$, $B$ and so on for these. By $a$, $b$ and
so on we mean a finite sequence of elements of $\M$.
By $a\in A$ we mean $a\in A^{length(a)}$.

By an automorphism we mean an automorphism of $\M$.
We write $Aut(A)$ for the set of all automorphisms of $\M$
such that $f\raj A=id_{A}$. By $S^{*}(A)$ we mean the the set of
all consistent complete types over $A$
and by $\tp (a,A)$ we mean the type of $a$ over $A$ in $\M$. $S^{m}(A)$
means the set
$\{ \tp (a,A)\vert\ a\in\M ,\ length(a)=m\}$ and
$S(A)=\cup_{m<\o}S^{m}(A)$.

We define $\k (\M )$
as $\k (T)$ is defined in the case of stable theories but for
strong splitting i.e. we let $\k (\M )$
be the least cardinal such that
there are no $a$, $b_{i}$ and $c_{i}$,
$i<\k (\M )$,
such that

(i) for all $i<\k (\M )$, there is an
infinite indiscernible set $I_{i}$ over
$\cup_{j<i}(b_{j}\cup c_{j})$ such that $b_{i},c_{i}\in I_{i}$,

(ii) for all $i<\k (\M )$, there is $\phi_{i}(x,y)$
such that
$\models\phi_{i}(a,b_{i})\wedge\neg\phi_{i}(a,c_{i})$.

We say that a type $p$ over $A$ is $\M$-consistent if there is $a\in\M$
such that $p\subseteq \tp (a,A)$ (i.e. there is $q\in S(A)$ such that
$p\subseteq q$).

\th 1.1 Lemma. ([Hy]) If $p\in S^{*}(A)$ is not
$\M$-consistent, then there is finite $B\subseteq A$
such that $p\raj B$ is not $\M$-consistent.

\th 1.2 Lemma.

(i) If $(a_{i})_{i<\o}$ is order-indiscernible over $A$ then
it is indiscernible over $A$.

(ii) Assume $\M$ is $\xi$-stable and $\vert I\vert >\xi\ge\vert A\vert$.
Then there is $J\subseteq I$ of power $>\xi$ such that it
is indiscernible over $A$.

(iii) If $I$ is infinite indiscernible over $A$ then for all
$\xi\le\vert\M\vert$
there is $J\supseteq I$ of power $\ge\xi$ such that $J$
is indiscernible over $A$.

(iv) For all indiscernible $I$ and $\phi (x,a)$, either
$X=\{ b\in I\vert\ \models\phi (b,a)\}$ or
$Y=\{ b\in I\vert\ \models\neg\phi (b,a)\}$ is of power $<\l (\M)$.

(v) There are no increasing sequence
of sets $A_{i}$, $i<\l (\M )$, and $a$
such that for all $i<\l (\M )$, $\tp (a,A_{i+1})$ splits over $A_{i}$.
So for all $A$ and $p\in S(A)$, there is $B\subseteq A$ of
power $<\l (\M )$, such that $p$ does not split over $B$.

(vi) For all $A$ and $p\in S(A)$, there is $B\subseteq A$ of
power $<\k (\M )$, such that $p$ does not split strongly over $B$.

\proof (i), (ii) and (v) as in [Hy]. (iii) follows immediately
from the homogeneity of $\M$. (vi) is trivial.

We prove (iv):
Assume not. Let $I$ be a counter example. Clearly we may assume that
$\vert I\vert =\l (\M )$.
Then By Lemma 1.1,
for
every $J\subseteq I$, the type
$$p_{J}=\{\phi (b,y)\vert\ b\in J\}\cup\{\neg\phi (b,y)\vert\ b\in I-J\}$$
is $\M$-consistent. Clearly
this contradicts $\l (\M )$-stability of $\M$.
$\eop$

\th 1.3 Corollary. $\k (\M )\le\l (\M )$.

\proof Follows immediately from Lemma 1.2 (v). $\eop$

We will
use Lascar strong types instead of strong types:

\th 1.4 Definition. Let $SE^{n}(A)$ be the set of all equivalence
relation $E$ in $\M^{n}$,
such that the number of equivalence classes is
$<\vert\M\vert$ and for all $f\in Aut(A)$, $a\ E\ b$ iff $f(a)\ E\ f(b)$.
Let $SE(A)=\cup_{n<\o}SE^{n}(A)$.

Notice that $E\in SE(A)$ need not be definable but
an indiscernible set over $A$ is also an indiscernible set for all
$E\in SE(A)$.

Usually we either do not mention the arities of the equivalence
relations we work with, or we mention that the arity is
f.ex. $m$, but we do not specify what $m$ is.
This is harmless since usually there is no danger of confusion.

\th 1.5 Lemma. If $I$ is an infinite indiscernible set over $A$,
then for all $E\in SE(A)$ and $a,b\in I$, $a\ E\ b$.

\proof Assume not. Let $E\in SE(A)$ be a counter example. Then
for all $a,b\in I$, $a\ne b$, $\neg (a\ E\ b)$. Then Lemma 1.2 (iii)
implies a contradiction with the number of equivalence classes of $E$.
$\eop$

\th 1.6. Lemma. If $E\in SE(A)$, $\vert A\vert\le\xi$ and $\M$ is
$\xi$-stable,
then the number of equivalence classes
of $E$ is $\le\xi$.

\proof Assume not. Then by Lemma 1.2 (ii), we can find
$I$ such that it is infinite indiscernible over $A$ and for all
$a,b\in I$, if $a\ne b$ then $\neg (a\ E\ b)$.
This contradicts Lemma 1.5. $\eop$

\th 1.7 Corollary. For all $A$ and $n<\o$,
there is $E_{min,A}^{n}\in SE^{n}(A)$ such that for all $a,b$ and
$E\in SE^{n}(A)$, $a\ E^{n}_{min,A}\ b$ implies
$a\ E\ b$.

\proof Clearly $\vert SE^{n}(A)\vert$ is restricted
($\le 2^{\vert S(A)\vert}$) and so
$\cap SE^{n}(A)\in SE(A)$. Trivially
$\cap SE^{n}(A)$ has the wanted property. $\eop$

\th 1.8 Definition.

(i) We say that $\A$ is $F^{\M}_{\k}$-saturated if for all
$A\subseteq\A$ of power $<\k$ and $a$, there is $b\in\A$ such
that $\tp (b,A)=\tp (a,A)$.

(ii) We say that $\A$ is strongly
$F^{\M}_{\k}$-saturated if for all $A\subseteq\A$ of power $<\k$
and $a$ of length $m$, there is $b\in\A$ such that
$b\ E\ a$ for all $E\in SE^{m}(A)$. We write $a$-saturated for strongly
$F^{\M}_{\k (\M )}$-saturated.

\th 1.9 Lemma.

(i) If $\A$ is strongly $F^{\M}_{\k}$-saturated then it is
$F^{\M}_{\k}$-saturated.

(ii) Assume $\vert A\vert\le\xi$, $\M$ is $\xi$-stable,
$\xi^{<\k}=\xi$ and there is a regular cardinal
$\d$ such that $\k\le\d\le\xi$.
Then there is strongly $F^{\M}_{\k}$-saturated
$\A\supseteq A$ such that $\vert\A\vert\le\xi$. Further more
if $\B\supseteq A$ is strongly $F^{\M}_{\k}$-saturated,
then we can choose $\A$ so that $\A\subseteq\B$.

(iii) Assume $\M$ is $\xi$-stable,
$\A$ is $F^{\M}_{\xi}$-saturated, $A\subseteq\A$
is of power $<\xi$ and $m<\o$. Then there are $a_{i}\in\A$, $i<\xi$, such that
for all $b$ of length $m$, there is
$i<\xi$ such that $a_{i}\ E\ b$, for all $E\in SE^{m}(A)$
i.e. $\A$ is strongly $F^{\M}_{\xi}$-saturated.

(iv) If $\A$ is $F^{\M}_{\l (\M )}$-saturated, then it is $a$-saturated.

(v) Assume
$\A$ is strongly
$F^{\M}_{\xi}$-saturated and $A\subseteq\A$
is of power $<\xi$. Then for all $B$ of power $<\xi$,
there is $f\in Aut(A)$ such that $f(B)\subseteq\A$ and for all
(finite sequences)
$b\in B$, $f(b)\ E^{m}_{min,A}\ b$.

\proof (i) is trivial.

(ii): For all $i\le\d$, choose sets $A_{i}$
of power $\le\xi$ as follows:
Let $A_{0}=A$ and if $i$ is limit then $A_{i}=\cup_{j<i}A_{j}$.
If $A_{i}$ is defined, then we let $A_{i+1}\supseteq A_{i}$
be such that for all $B\subseteq A_{i}$
of power $<\k$ and $a$ there is
$b\in A_{i+1}$ such that $b\ E^{m}_{min,B}\ a$. By Lemma 1.6,
we can find $A_{i+1}$ so that $\vert A_{i+1}\vert\le\xi$.
By Lemma 1.7, $A_{\d}$ is as wanted.

(iii): By Lemma 1.6, choose $b_{i}$, $i<\xi$, so that
for all $b$ there is $i<\xi$ such that $b\ E^{m}_{min,A}\ b_{i}$.
Since $\A$ is $F^{\M}_{\xi}$-saturated, we can choose $a_{i}\in\A$
so that there is $f\in Aut(A)$ such that for all $i<\xi$, $f(b_{i})=a_{i}$.
Clearly this implies the claim.

(iv): Immediate by (iii).

(v): For all $c\in B$, choose $a_{c}\in\A$ so that
$a_{c}\ E^{m}_{min,A}\ c$.
Since $\A$ is $F^{\M}_{\xi}$-saturated, there is
$f\in Aut(A\cup\{ a_{c}\vert\ c\in B\})$ such that $f(B)\subseteq\A$.
Clearly $f$ is as wanted.
$\eop$

\th 1.10 Definition. We write $f\in Saut(A)$ if $f\in Aut(A)$
and for all $a$, $f(a)\ E^{m}_{min,A}\ a$.

\th 1.11 Lemma. Assume $\M$ is $\xi$-stable and $\vert A\vert <\xi$. If
$a\ E^{m}_{min,A}\ b$, then there is $f\in Saut(A)$
such that $f(a)=b$.

\proof We define $a\ E\ b$ if there is $f\in Saut(A)$ such that $f(a)=b$.
Clearly it is enough to show that $E\in SE(A)$. For a contradiction,
assume that this is not the case. Since $E$ is an equivalence relation
and $f(E)=E$ for all $f\in Aut (A)$,
there are $a_{i}$, $i<\xi^{+}$, such that
for all $i\ne j$,
$\neg (a_{i}\ E\ a_{j})$. Choose $B\supseteq A$ of power
$\xi$ such that every $E^{m}_{min,A}$-equivalence class is
represented in $B$. Since $\M$ is $\xi$-stable, there are $i<j<\xi^{+}$,
such that $\tp (a_{i},B)=\tp (a_{j},B)$. Then there is $f\in Aut(B)$
such that $f(a_{i})=f(a_{j})$. By the choice of $B$, $f\in Saut(A)$,
a contradiction. $\eop$

\th 1.12 Lemma. Assume $\xi$ is such that for some $\xi'\ge\xi$, $\M$
is $\xi'$-stable.
If $\A$ is $F^{\M}_{\xi}$-saturated and
$A\subseteq\A$ has power
$<\xi$,
then $\tp (a,\A )$ does not split strongly over $A$ iff
for all $b,c\in\A$ and $\phi$, $b\ E^{m}_{min, A}\ c$ implies
$\models\phi (a,b)\leftrightarrow\phi (a,c)$.

\proof If $\tp (a,\A )$ splits strongly over $A$, then by Lemma 1.5,
there are $b,c\in\A$ and $\phi$, such that $b\ E^{m}_{min, A}\ c$ and
$\models\neg(\phi (a,b)\leftrightarrow\phi (a,c))$.
So we have proved the claim
from right to left.
We prove the other direction: For a contradiction assume that
there are $b,c\in\A$ and $\phi$, such that $b\ E^{m}_{min, A}\ c$ and
$\models\phi (a,b)\wedge\neg\phi (a,c)$.

We define an equivalence relation $E$
on $\M^{m}$ as follows:
$a\ E\ b$ if $a=b$ or there are $I_{i}$, $i<n<\o$,
such that they are infinite indiscernible over $A$, $a\in I_{0}$,
$b\in I_{n-1}$ and for all $i<n-1$, $I_{i}\cap I_{i+1}\ne\empty$.
Clearly $E$ is an equivalence relation and for all
$f\in Aut(A)$, $f(E)=E$. By Lemma 1.2 (ii), the number of
equivalence classes of $E$ is $<\vert\M\vert$.
So $E\in SE^{m}(A)$.

Then $b\ E\ c$ and $b\ne c$. Let $I_{i}$, $i<n$, be as in the definition
of $E$. Since $\A$ is $F^{\M}_{\vert A\vert^{+}+\o}$-saturated,
we may assume that for all $i<n$, $I_{i}\subseteq\A$.
Since $\tp (a,\A )$ does not split strongly over $A$, for all
$d\in I_{0}$, $\models\phi (a,d)$. So there is $d\in I_{1}$
such that $\models\phi (a,d)$. Again
since $\tp (a,\A )$ does not split strongly over $A$, for all
$d\in I_{1}$, $\models\phi (a,d)$. We can carry this on and finally
we get that $\models\phi (a,c)$, a contradiction. $\eop$

\th 1.13. Lemma. Assume $A\subseteq\A$, $\vert A\vert <\k (\M )$,
$\A$ is $a$-saturated and
$p\in S(\A )$ does not split strongly over $A$. Then for all
$B\supseteq\A$, there is $q\in S(B)$ such that $p\subseteq q$
and for all $C\supseteq B$ there is $r\in S(C)$,
which satisfies $q\subseteq r$ and
$r$ does not split strongly over $A$.

\proof  We define $q\in S^{*}(B)$ as follows:
$\phi (x,b)\in q$, $b\in B$, if there is $a\in\A$ such that
$a\ E^{m}_{min,A}\ b$ and $\phi (x,a)\in p$, where $m=length(b)$.
By Lemma 1.12, it is enough to show that $q$ is $\M$-consistent.
By Lemma 1.1, it is enough
to show that for all $a,a'\in\A$, if $a\ E^{m}_{min,A}\ a'$, then
$\phi (x,a)\in p$ implies $\phi (x,a')\in p$. This follows from
Lemma 1.12, since by Lemma 1.9 (i), $\A$ is
$F^{\M}_{\k (\M )}$-saturated.
$\eop$

\th 1.14 Lemma. Assume $A\subseteq\A\subseteq\B$,
$\vert A\vert <\k (\M )$,
$\B$ is $F^{\M}_{\k (\M )}$-saturated and for every $c\in\B$
there is $d\in\A$ such that $d\ E^{m}_{min,A}\ c$.
If $\tp (a,\A )=\tp (b,\A )$ and both $\tp (a,\B )$ and
$\tp (b,\B )$ do not split strongly over $A$, then
$\tp (a,\B )=\tp (b,\B )$.

\proof For a contradiction, assume $c\in\B$ and
$\models\phi (a,c)\wedge\neg\phi (b,c)$. Choose
$d\in\A$ such that $d\ E^{m}_{min,A}\ c$. By Lemma 1.12,
$\models\phi (a,d)\wedge\neg\phi (b,d)$, a contradiction. $\eop$

\th 1.15 Lemma. If
$\xi =\l (\M )+\xi^{<\k (\M )}$, then $\M$ is $\xi$-stable.

\proof Clearly we may assume that $\xi >\l (\M )$ and so
by Corollary 1.3, $\xi\ge\k (\M )^{+}$.
Let $A$ be a set of power
$\xi$. We show that $\vert S(A)\vert\le\xi$.

{\bf Claim.} There is $\A\supseteq A$ such that

(i) $\A$ is $F^{\M}_{\k (\M )}$-saturated,

(ii) $\vert\A\vert\le\xi$,

(iii) for all $B\subseteq\A$ of power $<\k (\M )$ there is
$\A_{B}\subseteq\A$
of power $\l (\M )$ satisfying: $B\subseteq\A_{B}$ and
for all $c\in\M$ there is
$d\in\A_{B}$ such that $d\ E^{m}_{min,A}\ c$.

\proof By induction on $i<\k (\M )^{+}$, we define $\A_{i}$
so that $\vert\A_{i}\vert\le\xi$, $A\subseteq\A_{0}$,
for $i<j$, $\A_{j}\subseteq\A_{i}$ and

(1) if $i$ is odd then for all $B\subseteq\cup_{j<i}\A_{j}$
of power $<\k (\M )$, there is
$\A_{B}\subseteq\A_{i}$
of power $\le\l (\M )$ satisfying: $B\subseteq\A_{B}$ and
for all $c\in\M$ there is
$d\in\A_{B}$ such that $d\ E^{m}_{min,A}\ c$,

(2) if $i$ is even then for all $B\subseteq\cup_{j<i}\A_{j}$
of power $<\k (\M )$, every $p\in S(B)$ is realized in $\A_{i}$.

\noindent
By Corollary 1.3, Lemma 1.6 and the fact that $\vert S(B)\vert\le\l (\M )$
for all $B$ of power $<\k (\M )^{+}$, it is easy to see that such
$\A_{i}$, $i<\k (\M )$, exist. Clearly $\A=\cup_{i<\k (\M )^{+}}\A_{i}$
is as wanted. $\eop$ Claim.

So it is enough to show that $\vert S(\A )\vert\le\xi$.
By Lemma 1.2 (vi),
for each $p\in S(\A )$, choose $B_{p}$ so that
$p$ does not split strongly over $B_{p}$ and
$\vert B_{p}\vert<\k (\M )$. Then
by Lemma 1.14,
every type $p\in S(\A )$
is determined by $p\raj \A_{B_{p}}$ and the fact that it
does not split strongly over $B$.
Since the number of possible $B$ is
$\xi^{<\k (\M )}=\xi$ and for each such $B$,
$\vert S(\A_{B})\vert\le\l (\M )$,
$\vert S(\A)\vert\le\xi\times\l (\M )=\xi$. $\eop$

\th 1.16 Lemma. If $\xi^{<\k (\M )}>\xi$, then $\M$ is not
$\xi$-stable.

\proof By the definition of $\l (\M )$, we may assume that $\xi\ge\l (\M )$.
Let $\k <\k (\M )$ be the least
cardinal such that $\xi^{\k}>\xi$. By the definition of $\k (\M )$,
there are  $a$, $b_{i}$ and $c_{i}$,
$i<\k$,
such that

(i) for all $i<\k$, there is an
infinite indiscernible set $I'_{i}$ over
$\cup_{j<i}(b_{j}\cup c_{j})$ such that $b_{i},c_{i}\in I'_{i}$,

(ii) for all $i<\k$, there is $\phi_{i}(x,y)$
such that
$\models\phi_{i}(a,b_{i})\wedge\neg\phi_{i}(a,c_{i})$.

{\bf Claim.} There are $I_{i}$, $i<\k$, such that
for all $i<\k$, $I_{i}=\{ d^{i}_{k}\vert\ k<\xi\}$
is indiscernible over $\cup_{j<i}I_{j}$,
$b_{i},c_{i}\in I_{i}$ and for $k<k'<\xi$, $d^{i}_{k}\ne d^{i}_{k'}$.

\proof By induction on $0<\a\le\k$, we define
$I^{\a}_{i}=\{ d^{\a ,i}_{k}\vert\ k<\xi\}$, $i<\a$,
such that

(1) for all $i<\a$, $I^{\a}_{i}$ is indiscernible over $\cup_{j<i}I^{\a}_{j}$
and $b_{i},c_{i}\in I^{\a}_{i}$,

(2) for all $\b <\a$, there is an automorphism $f$ such that
$f\raj\cup_{j<\b}(b_{j}\cup c_{j})=id_{\cup_{j<\b}(b_{j}\cup c_{j})}$
and for all $j<\b$, $f(d^{\b ,j}_{k})=d^{\a ,j}_{k}$, $k<\xi$,

(3) for all $i<\a$ and $k<k'<\xi$, $d^{\a ,i}_{k}\ne d^{\a ,i}_{k'}$.

\noindent
Clearly this is enough, since then $I^{\k}_{i}$, $i<\k$,
are as wanted.

By (2) and homogeneity of $\M$, limits are trivial, so we assume that
$\a =\b +1$ and that $I^{\b}_{j}$, $j<\b$, are defined.
By Lemma 1.15, there is
$\d >\xi$ such that $\M$ is $\d$-stable. By the assumptions
and Lemma 1.2 (iii), there is $J=\{ d_{k}\vert\ k<\d^{+}\}$
such that it is indiscernible over $\cup_{j<\b}(b_{j}\cup c_{j})$
and $b_{\b},c_{\b}\in J$. By Lemma 1.2 (ii), there is $I\subseteq J$
of power $\xi$, such that it is indiscernible over $\cup_{j<\b}I^{\b}_{j}$.
Since $J$ is indiscernible over $\cup_{j<\b}(b_{j}\cup c_{j})$,
there is an automorphism $f$ such that
$f\raj\cup_{j<\b}(b_{j}\cup c_{j})=id_{\cup_{j<\b}(b_{j}\cup c_{j})}$
and $b_{\b},c_{\b}\in\{ f(d)\vert\ d\in I\}$.
We let $I^{\a}_{\b}=f(I)$ and if $i<\b$, then $I^{\a}_{i}=f(I^{\b}_{i})$.
Clearly these are as required. $\eop$ Claim.

By Lemma 1.2 (iv) we may assume that for all $i<\k$,
$\models\phi_{i}(a,d^{i}_{k})$
iff $k=0$. Then for all $\n\in\xi^{\k}$ and $0<\a\le\k$, we define
function $f^{\n}_{\a}$ so that the following holds
($f^{\n}_{0}=id_{\M}$):

(a) for all $i<\b <\a$ and $\n\in\xi^{\k}$,
$f^{\n}_{\a}\raj I_{i}=f^{\n}_{\b}\raj I_{i}$,

(b) if $\a =\b +1$ and $\n\in\x^{\k}$,
then
$$f^{\n}_{\a}(f^{\n}_{\b}(d^{\b}_{0}))=f^{\n}_{\b}(d^{\b}_{\n (\b )}),$$
$$f^{\n}_{\a}(f^{\n}_{\b}(d^{\b}_{\n (\b )}))=f^{\n}_{\b}(d^{\b}_{0})$$
and for all $i<\xi$, $i\ne 0,\n (\b )$,
$$f^{\n}_{\a}(f^{\n}_{\b}(d^{\b}_{i}))=f^{\n}_{\b}(d^{\b}_{i}),$$

(c) if $\n\raj\a =\n'\raj\a$ then $f^{\n}_{\a}=f^{\n'}_{\a}$.

\noindent
It is easy to see that such $f^{\n}_{\a}$ exist. For limit $\a$
this follows from the homogeneity of $\M$ and for successors
this follows from the fact that
$f^{\n}_{\b}(I_{\b})$ is indiscernible over
$\cup_{i<\b}f^{\n}_{\b}(I_{i})$.

For all $\n\in\xi^{\k}$, let $a_{\n}=f^{\n}_{\k}(a)$.
Then clearly for $\n\ne\n '$, the types of $a_{\n}$ and
$a_{\n'}$ over
$A=\cup\{ f^{\nu}_{\a +1}(I_{\a})\vert\ \nu\in\xi^{\k},\ \a <\k\}$
are different. By the choice of $\k$, $\xi^{<\k}=\xi$ and so by (c),
$\vert A\vert =\xi$. Since $\xi^{\k}>\xi$, $\M$ is not
$\xi$-stable. $\eop$

So we have proved the following theorem. With slightly different
definitions this theorem is already proved in [Sh1].

\th 1.17 Theorem. $\M$ is $\xi$-stable iff
$\xi =\l (\M )+\xi^{<\k (\M )}$.

\proof Follows from Lemmas 1.15 and 1.16. $\eop$

Let
$\k_{r}(\M )$ be the least regular $\k\ge\k (\M )$.
By Lemma 1.16, $\l (\M )^{<\k (\M )}=\l (\M )$ and so
$cf(\l (\M ))\ge\k (\M )$. Because $cf(\l (\M ))$ is regular,
$\k_{r}(\M )\le\l (\M )$.

\chapter{2. Indiscernible sets}

In this chapter we prove basic properties of indiscernible sets.
We start by improving Lemma 1.2 (iv).

\th 2.1 Lemma. For all infinite indiscernible $I$ and $a$ there is
$p\in S(a)$ such that
$$\vert\{ b\in I\vert\ \tp (b,a)\ne p\}\vert <\k (\M ).$$

\proof Assume not. By Lemma 1.2 (iii), we may assume that
$I$ and $a$ are such that $I=\{ b_{i}\vert\ i<\k (\M )+\o\cdot\k (\M )\}$,
$b_{i}\ne b_{j}$ for $i\ne j$ and
for some $p\in S(a)$, $\tp (b_{i},a)=p$ iff $i\ge\k (\M )$.
For all $i<\k (\M )$, we define $A_{i}$ as follows:

(i) $A_{0}=\empty$,

(ii) $A_{i+1}=A_{i}\cup\{ b_{i-1}\}\cup
\{ b_{j}\vert\ \o\cdot i\le j<\o\cdot (i+1)\}$,

(iii) for limit $i$, $A_{i}=\cup_{j<i}A_{j}$.

\noindent
Then it is easy to see that for all $i<\k (\M )$
$\tp (a,A_{i+1})$ splits strongly over $A_{i}$, a contradiction. $\eop$

\th 2.2 Corollary. For all indiscernible $I$ and $\phi (x,a)$, either
$X=\{ b\in I\vert\ \models\phi (b,a)\}$ or
$Y=\{ b\in I\vert\ \models\neg\phi (b,a)\}$ is of power $<\k (\M)$.

\proof Follows immediately from Lemma 2.1.
$\eop$

\th 2.3 Definition. If $I$ is indiscernible and of power $\ge\k (\M )$,
we write $Av(I,A)$ for $\{\phi (x,a)\vert\ a\in A,
\ \phi\in L,\ \vert
\{ b\in I\vert\ \models\neg\phi (b,a)\}\vert <\k (\M )\}$.

\th 2.4 Lemma.

(i) If $I$ is indiscernible over $A$ and of power $\ge\k (\M )$,
then $I\cup\{ b\}$ is indiscernible over $A$ iff
$\tp (b,I\cup A)=Av(I,I\cup A)$.

(ii) If $I$ and $J$ are of power $\ge\k (\M )$ and $I\cup J$
is indiscernible, then for all $A$,
$Av(I,A)=Av(J,A)$.

(iii) If $I$ is indiscernible and of power $\ge\k (\M )$,
then for all $A$, $Av(I,A)$ is $\M$-consistent.

\proof (i) and (ii) are trivial. We prove (iii):
By (ii) and Lemma 1.2 (iii),
we may assume that
$\vert I\vert >\vert L\cup A\vert +\k_{r}(\M )$. Then
the claim follows by the pigeon hole principle
from (i). $\eop$

\th 2.5 Definition. Assume $I$ and $J$ are indiscernible sets of power
$\ge\k (\M )$.

(i) We say that $I$ is based on $A$ if for all
$B\supseteq A\cup I$, $Av(I,B)$ does not split strongly over $A$.

(ii) We say that $I$ and $J$ are equivalent if
for all $B$, $Av(I,B)=Av(J,B)$

(iii) We say that $I$ is stationary over $A$ if
$I$ is based on $A$ and for all
$f\in Aut(A)$, $f(I)$ and $I$ are equivalent.

\th 2.6 Lemma. Assume $I$ is an indiscernible set of power
$\ge\k (\M )$, $\vert A\vert <\xi$ and $\M$ is $\xi$-stable.
Then the following are equivalent:

(i) $I$ is based on $A$,

(ii) the number of non-equivalent indiscernible sets in
$\{ f(I)\vert\ f\in Aut(A)\}$ is $\le\xi$,

(iii) the number of non-equivalent indiscernible sets in
$\{ f(I)\vert\ f\in Aut(A)\}$ is $<\vert\M\vert$.

\proof (i)$\Rightarrow$(ii): Assume not. Let $f_{i}(I)$, $i<\xi^{+}$,
be a counter example.
For all $i<\l (\M )$, choose $\A_{i}$ so that

(a) $A\subseteq\A_{0}$ and
every type $p\in S(A)$ is realized
in $\A_{0}$,

(b) if $i<j$, then $\A_{i}\subseteq\A_{j}$ and
for limit $i$, $\A_{i}=\cup_{j<i}\A_{j}$,

(c) every type $p\in S(\A_{i})$ is realized
in $\A_{i+1}$,

(d) $\vert\A_{i}\vert\le\xi$.

\noindent
Let $\A =\cup_{i<\l (\M )}\A_{i}$.
Since $\M$ is
$\xi$-stable there are $i\ne j$ such that
$Av(f_{i}(I),\A )=Av(f_{j}(I),\A )$. Let $a$ be such that
$Av(f_{i}(I),\A\cup\{ a\} )\ne Av(f_{j}(I),\A\cup\{ a\} )$.
By Lemma 1.2 (v), choose $i<\l (\M )$
so that $\tp (a,\A_{i+\o})$ does not split over $\A_{i}$.
Without loss of generality, we may assume that $i=0$.
For all
$i<\o$, choose $a_{i}\in\A_{i+1}$ so that
$\tp (a_{i},\cup_{j\le i}\A_{j})=\tp (a,\cup_{j\le i}\A_{j})$.
By an easy induction, we see that $\{ a\}\cup\{ a_{i}\vert\ i<\o\}$ is
order-indiscernible over $\A$ and so also over $A$.
By Lemma 1.2 (i), $\{ a\}\cup\{ a_{i}\vert\ i<\o\}$ is
indiscernible over $A$. But then clearly either
$Av(f_{i}(I),\A\cup\{ a\} )$ or $Av(f_{j}(I),\A\cup\{ a\} )$
splits strongly over $A$, a contradiction.

(ii)$\Rightarrow$(iii): Trivial.

(iii)$\Rightarrow$(i): Assume not. Then by Lemma 1.2 (iii), we can find
$J=\{ a_{i}\vert\ i<\vert\M\vert\}$ and $\phi (x,y)$ such that
$J$ is indiscernible over $A$, for $i\ne j$, $a_{i}\ne a_{j}$, and
$\phi (x,a_{i})\in Av(I,J)$ iff $i=0$. But then for all
$i<\vert\M\vert$, we can find $f_{i}\in Aut(A)$ such that
for all $j<i$, $\phi (x,a_{j})\not\in Av(f_{i}(I),J)$ but
$\phi (x,a_{i})\in Av(f_{i}(I),J)$. Clearly these
$f_{i}(I)$ are not equivalent, a contradiction. $\eop$

\chapter{3. Independence}

In this chapter we define an independence relation and
prove the basic properties of it. This independence notion
is an improved version of the one defined in [Hy].
It satisfies weak versions of the basic properties of forking.
E.g. $a\da_{A}A$ holds assuming $A$ is $a$-saturated.

\th 3.1 Definition.

(i) We write $a\da_{A}B$ if
there is $C\subseteq A$ of power $<\k (\M )$ such that for all
$D\supseteq A\cup B$ there is $b$ which satisfies:
$\tp (b,A\cup B)=\tp (a,A\cup B)$ and
$\tp (b,D)$ does not split strongly over $C$. We write
$C\da_{A}B$ if for all $a\in C$, $a\da_{A}B$.

(ii) We say that $\tp (a,A)$ is bounded if
$\vert\{ b\vert\ \tp (b,A)=\tp (a,A)\}\vert <\vert\M\vert$.
If $\tp (a,A)$ is not bounded, we say that it is unbounded.

\th 3.2 Lemma.

(i) If $A\subseteq A'\subseteq B'\subseteq B$ and
$a\da_{A}B$ then $a\da_{A'}B'$.

(ii) If $A\subseteq B$ and $a\da_{A}B$ then for all
$C\supseteq B$ there is $b$ such that $\tp (b,B)=\tp (a,B)$
and $b\da_{A}C$.

(iii) Assume that $\A$ is $a$-saturated.
If $A\subseteq\A$ is such that
$\tp (a,\A )$ does not split strongly over $A$ then for all
$B$ such that $A\subseteq B\subseteq\A$,
$a\da_{B}\A$. Especially $a\da_{\A}\A$.

(iv) Assume $a$ and $A$ are such that
$\tp (a,A)$ is bounded.
Then for all $B\supseteq A$, $\tp (a,B)$
does not split strongly over $A$.

(v) Assume $A\subseteq B$
and $\tp (a,A)$ is unbounded. If
$\tp (a,B)$ is bounded,
then $a\nda_{A}B$.

(vi) Assume $\A$ is $a$-saturated
and $a\not\in\A$. Then
$\tp (a,\A )$ is unbounded.

(vii) Let $\xi$ be a cardinal. Assume $a$ and $A$ are such that
$\tp (a,A)$ is unbounded
and $a\da_{A}A$. If $a_{i}$, $i<\xi$, are such that
for all $i<\xi$, $\tp (a_{i},A)=\tp (a,A)$ and
$a_{i}\da_{A}\cup_{j<i}a_{j}$, then
$\vert\{ a_{i}\vert\ i<\xi\}\vert =\xi$.

(viii) Assume $A\subseteq B$, $a\da_{A}A$ and
$\tp (a,A)$ is unbounded.
Then there is $b$
such that $b\da_{A}B$ and $b\ E^{m}_{min,A}\ a$.

(ix) If $a\da_{A}b\cup c$ and $b\ E^{m}_{min,A}\ c$, then
$\tp (b,A\cup a)=\tp (c,A\cup a)$.

\proof (i) is immediate by the definition of $\da$.

(ii): Choose $a$-saturated
$\D\supseteq C$. Since $a\da_{A}B$, there are $b$ and
$A'\subseteq A$ such that $\tp (b,B)=\tp (a,B)$,
$\vert A'\vert<\k (\M )$ and
$\tp (b,\D )$ does not split strongly over $A'$. By Lemma 1.13,
$b$ is as wanted.

(iii): By Lemmas 1.2 (vi) and 1.13, $a\da_{A}\A$ and so by (i),
$a\da_{B}\A$.

(iv): Assume not. Then there are distinct $a_{i}$, $i<\vert\M\vert$,
and $\phi$, such that
$\{ a_{i}\vert\ i<\vert\M\vert\}$ is indiscernible over $A$ and
$\models\phi (a,a_{i})$ iff $i=0$. For all
$\k (\M )\le i<\vert\M\vert$, find an automorphism $f_{i}\in Aut(A)$
such that $f_{i}(a_{0})=a_{i}$, $f(a_{i})=a_{0}$ and for all
$0<j<i$, $f_{i}(a_{j})=a_{j}$. By Corollary 2.2, it is easy to see
that $\{ f_{i}(a)\vert\ \k (\M )\le i<\vert\M\vert\}$
contains $\vert\M\vert$ distinct elements, a contradiction.

(v): Assume not. Then by (ii)
we can find
$C\supseteq B$
and $b$ such that $\tp (b,B)=\tp (a,B)$, $b\da_{A}C$ and
$b\in C$. By Lemma 1.2 (ii), there is an infinite indiscernible
set $I$ over $A$ such that $b\in I$. Clearly we cannot find $c$
such that
$\tp (c,C)=\tp (b,C)$ and $\tp (c,C\cup I)$ does not split strongly over
some $A'\subseteq A$, a contradiction.

(vi): Follows immediately from (iii) and (v).

(vii): Immediate by (v).

(viii): Let $\xi>\vert A\vert$ be such that $\M$ is $\xi$-stable.
Choose $a_{i}$, $i<\xi^{+}$ so that $\tp (a_{i},A)=\tp (a,A)$
and $a_{i}\da_{A}\cup_{j<i}a_{j}$.
By (vii) and Lemma 1.2 (ii), we may assume that
$\{ a_{i}\vert\ i<\o\}$ is infinite indiscernible over $A$.
Clearly we may also assume that $a=a_{0}$. Let $d=a_{1}$. Then
$\tp (d,A)=\tp (a,A)$, $d\da_{A}a$ and by Lemma 1.5,
$d\ E^{m}_{min,A}\ a$. Then we can choose $b$ so that
$\tp (b,A\cup a)=\tp (d,A\cup a)$ and
$b\da_{A}a\cup B$. Clearly then $b$ is as wanted.

(ix) Follows immediately from Lemma 1.12. Notice that if
$b\ E^{m}_{min,A}\ c$, then for all $d\in A$,
$b\cup d\ E^{m+k}_{min,A}\ c\cup d$.
$\eop$

\th 3.3 Definition.

(i) We say that $\M$-consistent $p\in S(A)$ is stationary
if for all $a$, $b$ and $B\supseteq A$ the following holds:
if $\tp (a,A)=\tp (b,A)=p$, $a\da_{A}B$ and $b\da_{A}B$ then
$\tp (a,B)=\tp (b,B)$.

(ii) We say that $I$ is $A$-independent if for all
$a\in I$, $a\da_{A}I-\{ a\}$.

\th 3.4 Lemma. If $\A$ is $a$-saturated, then every
$\M$-consistent $p\in S(\A )$
is stationary.

\proof Assume not. Choose $\B\supseteq\A$, $a$ and $b$ so that
$\tp (a,\A )=\tp (b,\A )$, $a\da_{\A}\B$, $b\da_{\A}\B$ and $\tp (a,\B )\ne \tp (b,\B )$.
By Lemma 3.2 (ii) we may assume that $\B$ is
$F^{\M}_{\k (\M )}$-saturated. Choose
$c\in\B$ and $\phi$ so that $\models\phi (a,c)\wedge\neg\phi (b,c)$.
Let $A\subseteq\A$ be such that $\vert A\vert <\k (\M )$ and
both $\tp (a,\B )$ and $\tp (b,\B )$ do not split strongly over $A$.
Choose $d\in \A$ so that $d\ E^{m}_{min, A}\ c$. By Lemma 1.12,
a contradiction follows. $\eop$

\th 3.5 Corollary.

(i) Assume $\A$ is $a$-saturated. If $a\nda_{\A}B$, then there is $b\in B$
such that $a\nda_{\A}b$.

(ii) If $\A$ is $a$-saturated and $a_{i}$, $i<\a$, are such that
$a_{0}\not\in\A$,
for all $i,j$, $\tp (a_{i},\A )=\tp (a_{j},\A )$ and $a_{i}\da_{\A}\cup_{j<i}a_{j}$,
then $\{ a_{i}\vert\ i<\a\}$ is indiscernible over $\A$ and
$\A$-independent and if
$i\ne j$, then $a_{i}\ne a_{j}$.

(iii) Assume $\A$ is $a$-saturated. Then for all $B\supseteq\A$ and
$C$ there is $D$ such that $\tp (D,\A )=\tp (C,\A )$ and
$D\da_{\A}B$.

(iv) If $A\subseteq\B\subseteq C$, $\B$ is $a$-saturated,
$a\da_{A}\B$ and $a\da_{\B}C$, then $a\da_{A}C$.

(v) Assume $\A$ is $a$-saturated, $\tp (a,\A )$ does not split strongly over
$A\subseteq\A$ and $\vert A\vert <\k (\M )$. Then $a\nda_{\A}B$ iff
there is finite $b\in \A\cup B$ such that $a\nda_{A}b$.

\proof (i) follows immediately from Lemma 3.4 (if $a\nda_{\A}B$,
then $\tp (a,\A\cup B)$ is not the unique free extension of
$\tp (a,\A)$, which can be detected from a finite sequence).

(ii): By Lemma 3.4, it is easy to see that
$\{ a_{i}\vert\ i<\a\}$ is order-indiscernible over $\A$.
By Lemma 1.2 (i),
$\{ a_{i}\vert\ i<\a\}$ is indiscernible over $\A$.
Clearly this implies that $\{ a_{i}\vert\ i<\a\}$ is
$\A$-independent. The last claim follows from Lemma 3.2 (v).

(iii): Clearly it is enough to prove the following:
If $D\da_{\A}B$, then for all $c$ there is $d$ such that
$\tp (d,\A\cup D)=\tp (c,\A\cup D)$ and $d\cup D\da_{\A}B$.
This follows from Lemmas 1.1, 3.2 (ii) and 3.4.

(iv): Choose $b$ so that
$\tp (b,\B )=\tp (a,\B )$ and $b\da_{A}C$. Then $b\da_{\B}C$ and so
by Lemma 3.4,
we get $\tp (b,C)=\tp (a,C)$. Clearly this implies the claim.

(v): If $a\da_{\A}B$ then by (iv), $a\da_{A}\A\cup B$ from which
it follows that there are no finite $b\in \A\cup B$ such that $a\nda_{A}b$.
On the other hand if $a\nda_{\A}B$, then $\tp (a,\A\cup B)$
is not the unique 'free' extension of $\tp (a,\A )$ defined
in the proof of Lemma 1.13. This means that there are
$c\in\A$ and $d\in\A\cup B$ such that $c\ E^{m}_{min,A}\ d$
and $\tp (c,A\cup a)\ne \tp (d,A\cup a)$. Clearly $a\nda_{A}c\cup d$.
$\eop$

\th 3.6 Lemma. If $\A$ is $a$-saturated and $a\da_{\A}b$, then
$b\da_{\A}a$.

\proof Assume not. Let $\xi >\vert\A\vert$ be such that
$\M$ is $\xi$-stable. For all $i<\xi^{+}$, choose
$a_{i}$ and $b_{i}$ so that $\tp (a_{i},\A )=\tp (a,\A )$,
$a_{i}\da_{\A}\cup_{j<i}(a_{j}\cup b_{j})$,
$\tp (b_{i},\A )=\tp (b,\A )$ and
$b_{i}\da_{\A}a_{i}\cup\bigcup_{j<i}(a_{j}\cup b_{j})$. Then by
Lemma 3.4, $b_{i}\nda_{\A}a_{j}$ iff $j>i$. Clearly
this contradicts Lemma 1.2 (ii). $\eop$

\th 3.7 Corollary. For all $a,b$ and $A$, $b\da_{A}A$ and
$a\da_{A}b$ implies
$b\da_{A}a$.

\proof Assume not. Choose $a$-saturated $\A\supseteq A$ and
$b'$ so that $\tp (b',A)=\tp (b,A)$ and $b'\da_{A}\A$. We may assume that $b'=b$.
Then choose $a'$ so that $\tp (a'.A\cup b)=\tp (a,A\cup b)$ and
$a'\da_{A}\A\cup b$. By Lemma 3.6, $b\da_{\A}a'$. By Corollary 3.5 (iii),
$b\da_{A}a'$ and so $b\da_{A}a$. $\eop$

\th 3.8 Lemma.

(i) If $b\da_{A}D$ and $c\da_{A\cup b}D$, then $b\cup c\da_{A}D$.

(ii) If $\A$ is $a$-saturated, $B\da_{\A}D$ and $C\da_{\A\cup B}D$,
then $B\cup C\da_{\A}D$.

(iii) Assume $\A$ is $a$-saturated and $B\supseteq\A$.
If $a\da_{\A}B$, $a\da_{B}C$
and there is $D\subseteq B$ (f.ex. $D=B$) such that $C\da_{D}B$,
then $a\da_{\A}B\cup C$.

(iv) Assume $\A$ is $a$-saturated. If $a\da_{\A}b$ and
$a\cup b\da_{\A}B$, then $a\da_{\A}B\cup b$.

(v) Assume $a\da_{A}A$, for all $i<\o$, $\tp (a_{i},A)=\tp (a,A)$
and $a_{i}\da_{A}\cup_{j<i}a_{j}$. Then for all $n<\o$,
$\{ a_{i}\vert\ i<n\}$ is $A$-independent.

\proof (i): Choose $B\subseteq A$ of power $<\k (\M )$
such that

(a) for all
$C\supseteq A\cup D$ there is $b'$ which satisfies:
$\tp (b',A\cup D)=\tp (b,A\cup D)$ and
$\tp (b',C)$ does not split strongly over $B$

\noindent
and

(b) for all $C\supseteq A\cup D\cup b$ there is $c'$ which satisfies:
$\tp (c',A\cup D\cup b)=\tp (c,A\cup D\cup b)$ and
$\tp (c',C)$ does not split strongly over $B\cup b$.

\noindent
Let $C\supseteq A\cup D$ be arbitrary. Choose $b'$ as in (a) above.
By (b) above we can find $c'$ such that
$\tp (c'\cup b',A\cup D)=\tp (c\cup b,A\cup D)$ and
$\tp (c',C\cup b')$ does not split strongly over $B\cup b'$.

For a contradiction, assume $\tp (b'\cup c',C)$ splits
strongly over $B$. Let $I=\{ a_{i}\vert\ i<\o\}\subseteq C$ and $\phi$
be such that $I$ is
indiscernible over $B$ and

(c) $\models\phi (c',b',a_{0})\wedge\neg\phi (c',b',a_{1})$.

{\bf Claim.} $I$ is indiscernible over $B\cup b'$.

\proof If not, then (change the enumeration if necessary)
there is $\psi$ over $B$ such that
$\models\psi (b',a_{0},...,a_{n-1})\wedge
\neg\psi (b',a_{n},...,a_{2n-1})$.
Since
$$\{ (a_{m\cdot n},...,a_{(m+1)\cdot n-1})\vert\ m<\o\}$$
is indiscernible over $B$, we have a contradiction
with the choice of $b'$. $\eop$ Claim.

By Claim and (c), $\tp (c',C\cup b')$ splits strongly over $B\cup b'$.
This contradicts the choice of $c'$.

(ii): Clearly we may assume that $C$ is finite.
Let $b\in B$ be arbitrary. We show that $C\cup b\da_{\A}D$.
Choose $A\subseteq\A$ and $A'\subseteq B$ such that

(a) $b\in A'$, $\vert A\cup A'\vert <\k (\M )$,

(b) for all
$D'\supseteq\A\cup B\cup D$ there is $C'$ which satisfies:
$\tp (C',\A\cup B\cup D)=\tp (C,A\cup B\cup D)$ and
$\tp (C',D')$ does not split strongly over $A\cup A'$

(c) for all
$D'\supseteq\A\cup D$ and $a\in A'$,
there is $a'$ which satisfies:
$\tp (a',\A\cup D)=\tp (a,A\cup D)$ and
$\tp (a',D')$ does not split strongly over $A$.

\noindent
Then we can proceed as in (i). (We assume that $\A$ is $a$-saturated
in order to be able to use Corollary 3.5 (iii).)

(iii): By Lemma 3.6, $B\da_{\A}a$. By Corollary 3.7,
$C\da_{B}a$. By (ii), these imply $B\cup C\da_{\A}a$,
from which we get the claim by Lemma 3.6.

(iv): Choose $a'$ so that $\tp (a',\A\cup b)=\tp (a,\A\cup b)$ and
$a'\da_{\A}B\cup b$. By (i) and Lemma 3.4,
$\tp (a'\cup b,\A\cup B)=\tp (a\cup b,\A\cup B)$.

(v): By (i) it is easy to see that

(*) for all $n<\o$,
$\cup_{i<n}a_{i}\da_{A}A$.

\noindent
We prove the claim by induction on $n$.
For $n=1$ the claim follows immediately from the assumptions.
Let $i<n$. We show that $a_{i}\da_{A}\cup\{ a_{j}\vert\ j<n,\ j\ne i\}$.
If $i=n-1$, then this is assumption. So assume that $i<n-1$.
By the choice of $a_{n-1}$,
$$a_{n-1}\da_{A\cup\bigcup\{ a_{j}\vert\ j<n-1,\ j\ne i\}}a_{i}.$$
By the induction assumption
$$a_{i}\da_{A}\cup\{ a_{j}\vert\ j<n-1,\ j\ne i\}$$
and by (*) and Corollary 3.7
$$\cup\{ a_{j}\vert\ j<n-1,\ j\ne i\}\da_{A}a_{i}.$$
By (i),
$$a_{n-1}\cup\bigcup\{ a_{j}\vert\ j<n-1,\ j\ne i\}\da_{A}a_{i}.$$
By Corollary 3.7, the claim follows.
$\eop$

\th 3.9 Lemma. Assume
$B\supseteq A$ and
$\tp (a,A)$ is unbounded.
Then $a\da_{A}B$ iff there is
an indiscernible set $I$ over $A$ such that
$\vert I\vert\ge\k (\M )$, $I$ is based on
some $A'\subseteq A$ of power $<\k (\M )$ and
$Av(I,B)=\tp (a,B)$.

\proof From right to left the claim is trivial. So we prove the
other direction. Without loss of generality, we may assume that
$B$ is $a$-saturated.
Let $A'\subseteq A$ be such that $\vert A'\vert <\k (\M )$
and for all $C\supseteq B$ there is $b$ such that
$\tp (b,B)=\tp (a,B)$ and $\tp (b,C)$ does not split strongly over $A'$.
Let $\xi >\vert B\vert$ be a regular cardinal such that
$\M$ is $\xi$-stable. For all $i<\xi^{+}$ we define
$\B_{i}$ and $a_{i}$ so that

(i) $\B_{i}$, $i<\xi^{+}$, is an increasing sequence of
$\xi$-saturated models of power $\xi$ and $B\subseteq\B_{0}$,

(ii) for all $i<\xi^{+}$, $\tp (a_{i},B)=\tp (a,B)$, $a_{i}\in\B_{i+1}-\B_{i}$
and $\tp (a_{i},\B_{i})$ does not split strongly over $A'$
(so $a_{i}\da_{A'}\B_{i}$).

\noindent
By Lemma 3.2 (v) and
Corollary 3.5 (ii), $\{ a_{i}\vert\ i<\xi^{+}\}$
is indiscernible over $B$ and $a_{j}\ne a_{j}$
for all $i<j<\xi^{+}$. We prove that
$I=\{ a_{i}\vert\ i<\k (\M )\}$ is as wanted.

Clearly it is enough to show that $I$ is based on $A'$. For a contradiction,
assume that $C\supseteq B$ is such that $Av(I,C)$ splits strongly over
$A'$. Clearly we may assume that $C\subseteq\B_{\k (\M )+1}$.
By Lemma 1.2 (ii) there is $J\subseteq\xi^{+}-(\k (\M )+1)$, such that
$\vert J\vert =\xi^{+}$ and
$\{ a_{i}\vert\ i\in J\}$ is indiscernible over $C$.
Then $\tp (a_{i},C)=Av(I,C)$ for all $i\in J$. By (ii) above,
for all $i\in J$,
$\tp (a_{i},C)$ does not split strongly over $A'$, a contradiction. $\eop$

\th 3.10 Lemma. Assume $a\ E^{m}_{min,A}\ b$, $a\da_{A}c$ and
$b\da_{A}c$. If $c\da_{A}A$ or
$\tp (a,A)$ is bounded or
$\tp (c,A)$ is bounded,
then $\tp (a,A\cup c)=\tp (b,A\cup c)$.

\proof We divide the proof to three cases:

Case 1. $\tp (c,A)$ is bounded:
Let $B$ be the set of all $e$ such that
$\tp (e,A)$ is bounded.
Then $\vert B\vert<\vert\M\vert$ and so
$\vert S(A\cup B)\vert <\vert\M\vert$. We define $E$ so that
$x\ E\ y$ if $\tp (x,A\cup B)=\tp (y,A\cup B)$.
Since for all $f\in Aut(A)$, $f(A\cup B)=A\cup B$, $E\in SE(A)$.
Clearly this implies the claim.

Case 2. $\tp (a,A)$ is bounded:
Define $E$ so that $x\ E\ y$ if $x=y$ or $\tp (x,A)\ne \tp (a,A)$ and
$\tp (y,A)\ne \tp (a,A)$. Clearly $E\in SE^{m}(A)$, and so $a=b$ from
which the claim follows.

Case 3. $\tp (a,A)$ is unbounded and $c\da_{A}A$:
Assume the claim is not true.
Let $\xi>\vert A\vert$ be such that $\M$ is $\xi$-stable.
Choose $a_{i}$, $i<\xi^{+}$ so that $\tp (a_{i},A\cup c)=\tp (a,A\cup c)$
and $a_{i}\da_{A}c\cup\bigcup_{j<i}a_{j}$.
By Lemmas 3.2 (vii) and 1.2 (ii), we may assume that
$\{ a_{i}\vert\ i<\o\}$ is infinite indiscernible over $A$.
Clearly we may also assume that $a=a_{0}$. Let $d=a_{1}$. Then
$\tp (d,A\cup c)=\tp (a,A\cup c)$, $d\da_{A}a\cup c$ and by Lemma 1.5,
$d\ E^{m}_{min,A}\ a$. Then we can choose this $d$ so that
in addition, $d\da_{A}a\cup c\cup b$. By Lemma 3.8 (i),
$b\cup d\da_{A}c$. By Corollary 3.7, $c\da_{A}b\cup d$.
Since $d\ E^{m}_{min,A}\ b$,
this contradicts Lemma 3.2 (ix). $\eop$

Notice that in the case(s) 1 (and 2) above the assumptions
$a\da_{A}c$ and $b\da_{A}c$ are not used.

\th 3.11 Corollary. Assume $a_{i}$, $i<\o$, are such that
for all $i,j<\o$, $a_{i}\ E^{m}_{min,A}\ a_{j}$
and for all $i<\o$, $a_{i}\da_{A}\cup_{j<i}a_{j}$.
Then for all $i\ne j$, $a_{i}\ne a_{j}$ and
$\{ a_{i}\vert\ i<\o\}$ is indiscernible over $A$.

\proof By Lemma 3.2 (vii), for all $i\ne j$, $a_{i}\ne a_{j}$.
We show that for all $i_{0}<i_{1}<...<i_{n}<\o$,
$\tp (a_{0}\cup ...\cup a_{n},A)=
\tp (a_{i_{0}}\cup ...\cup a_{i_{n}},A)$. By Lemma 1.2 (i),
this is enough.

By Lemma 3.8 (v),
$\{ a_{i}\vert\ i\le i_{n}\}$ is $A$-independent and
by Lemma 3.8 (i), it is easy to see that
$\cup\{ a_{i}\vert\ i\le i_{n}\}\da_{A}A$.
So by Lemma 3.10,
$\tp (a_{0},A\cup\bigcup_{0<k\le n}a_{i_{k}})=
\tp (a_{i_{0}},A\cup\bigcup_{0<k\le n}a_{i_{k}})$.
So it is enough to show that
$\tp (a_{0}\cup ...\cup a_{n},A)=
\tp (a_{0}\cup a_{i_{1}}\cup ...\cup a_{i_{n}},A)$.
As above we can see that
$\tp (a_{1},A\cup a_{0}\cup\bigcup_{1<k\le n}a_{i_{k}})=
\tp (a_{i_{1}},A\cup a_{0}\cup\bigcup_{1<k\le n}a_{i_{k}})$.
So it is enough to show that
$\tp (a_{0}\cup ...\cup a_{n},A)=
\tp (a_{0}\cup a_{1}\cup a_{i_{2}}\cup ...\cup a_{i_{n}},A)$.
We can carry this on and get the claim. $\eop$

\th 3.12 Theorem. Assume $a\da_{A}c$, $b\da_{A}c$ and
$a\ E^{m}_{min,A}\ b$. Then $\tp (a,A\cup c)=\tp (b,A\cup c)$.

\proof Assume not. As in the proof of Lemma 3.10 (Case 3.),
we can find $a'$ and $b'$ such that
$\tp (a',A\cup c)=\tp (a,A\cup c)$, $\tp (b',A\cup c)=\tp (b,A\cup c)$,
$a'\da_{A}c\cup a$, $b'\da_{A}c\cup b$,
$a'\ E^{m}_{min,A}\ a$ and $b'\ E^{m}_{min,A}\ b$.
For all $i<\k (\M )$, choose $a_{i}$ so that
$a_{i}\da_{A}c\cup a\cup b\cup\bigcup_{j<i}a_{j}$,
if $i$ is odd, then $\tp (a_{i},A\cup c\cup a)=\tp (a',A\cup c\cup a)$
and if $i$ is even, then $\tp (a_{i},A\cup c\cup b)=\tp (b',A\cup c\cup b)$.
By Corollary 3.11, for all $i\ne j$, $a_{i}\ne a_{j}$ and
$\{ a_{i}\vert\ i<\k (\M )\}$ is indiscernible over $A$.
Clearly this contradicts Lemma 2.1. $\eop$

\th 3.13 Lemma. Assume $\M$ is $\xi$-stable and $\vert A\vert\le\xi$.
Then there is $a$-saturated
$\A\supseteq A$ of power $\le\xi$.

\proof Immediate by Lemma 1.9 (ii) and the fact that
$\k_{r}(\M )\le\l (\M )$ is regular. $\eop$

\th 3.14 Theorem. Assume $\M$ is $\xi$-stable and $\vert A\vert\le\xi$.
Then there is $F^{\M}_{\xi}$-saturated
$\A\supseteq A$ of power $\le\xi$.

\proof By Lemma 3.13, there is an increasing
continuous sequence $A_{i}$, $i\le\xi\cdot\xi$,
of models of power $\le\xi$ such that

(i) $A\subseteq A_{0}$ and for all $i\le\xi\cdot\xi$,
$A_{i+1}$ is $a$-saturated,

(ii) for all $i<\xi\cdot\xi$ and $a$, there is $b\in A_{i+1}$
such that $\tp (b,A_{i})=\tp (a,A_{i})$.

\noindent
We show that $\A =A_{\xi\cdot\xi}$ is as wanted.
For this let $B\subseteq\A$ of power $<\xi$ and $b$
be arbitrary. We show that $\tp (b,B)$ is realized in $\A$.

By Theorem 1.17, $cf (\xi )\ge\k_{r}(\M )$ and so
$\A$ is $a$-saturated and there is
$\a' <\xi$ such that $b\da_{A_{\xi\cdot\a'}}\A$.
By the pigeon hole principle there is $\a <\xi$ such that
$\a\ge\a'$ and
$(A_{\xi\cdot (\a +1)}-A_{\xi\cdot\a})\cap B=\empty$.

{\bf Claim.} There is $\b <\xi$ such that
$B\da_{A_{\xi\cdot\a+\b}}A_{\xi\cdot\a+\b+1}$.

\proof Assume not. Then by the pigeon hole principle,
we can find $c\in B$ such that
$$\vert\{ \g<\xi\vert\ c\nda_{A_{\xi\cdot\a+\g}}A_{\xi\cdot\a+\g+1}\}
\vert\ge cf(\xi ).$$
But this is impossible
by Lemma 3.2 (iii), because
$cf(\xi)\ge\k_{r}(\M )$
and $\A_{\xi\cdot \g}$ is $a$-saturated for
all $\g\le\xi$ such that $cf(\g )\ge\k_{r}(\M )$. $\eop$ Claim.

Choose $c\in A_{\xi\cdot\a+\b +1}$ so that
$\tp (c,A_{\xi\cdot\a+\b})=\tp (b,A_{\xi\cdot\a+\b})$.
By Claim, $B\da_{A_{\xi\cdot\a+\b}}c$ and so
$c\da_{A_{\xi\cdot\a+\b}}B$. Since $b\da_{A_{\xi\cdot\a+\b}}B$,
Lemma 3.4 implies,
$\tp (c,A_{\xi\cdot\a+\b}\cup B)=\tp (b,A_{\xi\cdot\a+\b}\cup B)$. $\eop$

We finish this chapter by proving that
over $F^{\M}_{\l (\M )}$-saturated models
our independence notion is equivalent with
the notion used in [Hy].

\th 3.15 Lemma. Assume $\A$ is $F^{\M}_{\l (\M )}$-saturated model
and $B\supseteq\A$.
Then the following are equivalent:

(i) $a\da_{\A}B$.

(ii) For all $b\in B$ there is $A\subseteq\A$ of power $<\l (\M )$
such that $\tp (a,\A\cup b)$ does not split over $A$.

\proof Let $p\in S(\A )$ be arbitrary $\M$-consistent type.
Let $a$ be such that $\tp (a,\A )=p$ and $a\da_{\A}B$.
Let $a'$ be such that $\tp (a',\A )=p$ and
for all $b\in B$ there is $A\subseteq\A$ of power $<\l (\M )$
such that $\tp (a',\A\cup b)$ does not split over $A$.
We show that then $\tp (a,B)=\tp (a',B)$. This implies the claim, since
for all $\M$-consistent $p\in S(\A )$ such $a$ and $a'$
exist: The existence of $a$ follows from Lemma 3.2 (ii) and (iii)
and the existence of $a'$ can be seen as in [Hy].

For a contradiction, assume that there is $b\in B$ such that
$\tp (a,\A\cup b)\ne \tp (a',\A\cup b)$. By the choice of $a$ and $a'$
and Lemma 1.2 (vi), there is $A\subseteq\A$ of power $<\l (\M )$
such that $\tp (a,\A\cup b)$ does not split strongly over $A$,
$\tp (a',\A\cup b)$ and $\tp (b,\A )$ do not split over $A$
and $\tp (a,A\cup b)\ne \tp (a',A\cup b)$.
For all $i<\o$, choose $b_{i}\in\A$ so that
$\tp (b_{i},A\cup\bigcup_{j<i}b_{j})=\tp (b,A\cup\bigcup_{j<i}b_{j})$.
Since $\tp (b,\A )$ does not split over $A$,
by Lemma 1.2 (i), it is easy to see that
$\{ b_{i}\vert\ i<\o\}\cup\{ b\}$ is infinite indiscernible over $A$.
Since $\tp (a,\A )=\tp (a',\A )$, either $\tp (a,\A\cup b)$ or $\tp (a',\A\cup b)$
splits strongly over $A$, a contradiction. $\eop$

\chapter{4. Orthogonality}

In this chapter we study orthogonality. Since we do not have full
transitivity of $\da$, we need stationary pairs:

\th 4.1 Definition. Assume $A\subseteq B$ and $p\in S(B)$.
We say that $(p,A)$ is stationary pair if
for all $a$, $\tp (a,B)=p$ implies $a\da_{A}B$ and
for all $C\supseteq B$,
$a$ and $b$, the following holds: if $a\da_{A}C$, $b\da_{A}C$ and
$\tp (a,B)=\tp (b,B)=p$, then $\tp (a,C)=\tp (b,C)$.

\th 4.2 Lemma.

(i) Assume $A\subseteq B\subseteq C$, $a\da_{A}C$ and
$(\tp (a,B),A)$ is a stationary pair. Then $(\tp (a,C),A)$ is
a stationary pair.

(ii) Assume $A\subseteq B\subseteq C\subseteq D$,
$a\da_{A}C$, $a\da_{B}D$ and $(\tp (a,C),B)$ is a stationary pair.
Then $a\da_{A}D$.

\proof (i) is trivial, so we prove (ii): Choose $a'$
so that $\tp (a',C)=\tp (a,C)$ and $a'\da_{A}D$. Then
$a'\da_{B}D$ and so $\tp (a',D)=\tp (a,D)$ from which the claim follows.
$\eop$

\th 4.3 Lemma. Assume $\A$ is
$a$-saturated, $\tp (a,\A )$ does not split strongly over $A\subseteq\A$
and $\vert A\vert <\k (\M )$.
Then
there is $B\subseteq\A$
such that $A\subseteq B$,
$\vert B-A\vert <\o$, $B\da_{A}A$ and
$(\tp (a,B),A)$ is a stationary pair.

\proof By Lemma 1.13, $a\da_{A}\A$.
Choose $b_{i}$, $i\le\o$, so that for all $i\le\o$,
$\tp (b_{i},\A )=\tp (a,\A )$ and $b_{i}\da_{A}\A\cup\bigcup_{j<i}b_{j}$.
Then $\{ b_{i}\vert\ i\le\o\}$
is indiscernible over $A$ and by Lemma 3.8 (ii),
$$(*)\ \ \ \ \ \{ b_{i}\vert\ i<\o\}\da_{A}\A.$$
Especially,
$$(**)\ \ \ \ \ \{ b_{i}\vert\ i<\o\}\da_{A}A.$$
Without loss of generality, we may assume that $b_{\o}=a$.
Choose $a^{*}\in\A$
so that
$a^{*}\ E^{m}_{min,A}\ a$.
Let $B=A\cup a^{*}$
and $I=\{ b_{i}\vert\ i<\o\}$. Then $B\da_{A}A$.

{\bf Claim.} Assume $J\supseteq I$ is indiscernible over $A$,
$\tp (b,B)=\tp (a,B)$ and $b\da_{A}B\cup J\cup a$.
Then $J\cup\{ b\}$ is indiscernible over $A$.

\proof By Lemmas 1.12 and 1.5
it is enough to show that $\tp (b,A\cup I)=\tp (a,A\cup I)$.
By $(*)$,
$I\da_{A}a^{*}$. By the choice of $a^{*}$,
$a^{*}\da_{A}A$ and so by Corollary 3.7,
$a^{*}\da_{A}I$. By the choice of $b$ and Lemma 3.2 (i),
$b\da_{A\cup a^{*}}I$. By Lemma 3.8 (i),
$b\cup a^{*}\da_{A}I$. By
$(**)$ and Corollary 3.7,
$I\da_{A}a^{*}\cup b$. So by Lemma 3.2 (ix),
$\tp (a^{*},A\cup I)=\tp (b,A\cup I)$. Similarly we can see that
$I\da_{A}a^{*}\cup a$ and so by Lemma 3.2 (ix),
$\tp (a^{*},A\cup I)=\tp (a,A\cup I)$.
$\eop$ Claim.

We show that $(\tp (a,B),A)$
is a stationary pair. Assume not.
Since $\A$ is $F^{\M}_{\k (\M )}$-saturated,
we can find $b$ such that
$b\da_{A}\A$, $\tp (b,B)=\tp (a,B)$ and $\tp (b,\A )\ne \tp (a,\A )$.
Choose $c_{i}$, $i<\k (\M )$, so that for all $i<\k (\M )$,
$\tp (c_{i},\A )=\tp (b,\A )$ if $i$ is odd, $\tp (c_{i},\A )=\tp (a,\A )$ if
$i$ is even and for all $i<\k (\M )$,
$c_{i}\da_{A}\A\cup I\cup\bigcup_{j<i}c_{j}$.
By Claim $\{ c_{i}\vert\ i<\k (\M )\}$
is indiscernible. This contradicts Corollary 2.2. $\eop$

\th 4.4 Definition.

(i) We say that $p\in S(A)$ is orthogonal to
$q\in S(C)$ if for all $a$-saturated
$\A\supseteq A\cup C$ the following holds: if
$\tp (b,C)=q$, $b\da_{C}\A$,
$\tp (a,A)=p$ and $a\da_{A}\A$, then $a\da_{\A}b$.
We say that $p\in S(A)$ is orthogonal to $C$ if it is
orthogonal to every $q\in S(C)$.

(ii) We say that a stationary pair $(p,A)$ is
orthogonal to $q\in S(C)$ if for all $a$-saturated
$\A\supseteq C\cup dom(p)$
the following holds: if $\tp (b,C)=q$,
$b\da_{C}\A$, $\tp (a,dom(p))=p$ and
$a\da_{A}\A$, then $a\da_{\A}b$.
We say that a stationary pair $(p,A)$ is
orthogonal to $C$ if it is orthogonal to every $q\in S(C)$.

\th 4.5 Lemma. Assume $\A$ is $a$-saturated, $A\subseteq B\subseteq\A$,
$a\da_{A}\A$ and $(\tp (a,B),A)$ is a stationary pair. Then
$\tp (a,\A )$ is orthogonal to $C$ iff $(\tp (a,B),A)$ is orthogonal to $C$.

\proof Immediate. $\eop$

\th 4.6 Lemma. Assume $A\subseteq\A$, $\A$ is $a$-saturated
and $p\in S(\A )$. Then the following are equivalent.

(i) $p$ is orthogonal to $A$.

(ii) For all $a$ and $b$, if $\tp (a,\A )=p$ and $b\da_{A}\A$, then
$a\da_{\A}b$.

\proof Clearly (i) implies (ii) and so we prove the other direction.
Assume (ii) and for a contradiction assume that there is
$a$-saturated $\C\supseteq\A$ and $a$ and $b$ such that
$\tp (a,\A )=p$, $a\da_{\A}\C$, $b\da_{A}\C$
and $a\nda_{\C}b$.

Choose $B_{0}\subseteq B_{1}\subseteq\A$ so that

(1) $\vert B_{1}\vert <\k (\M )$,

(2) $a\da_{B_{0}}\A$ and $b\da_{B_{0}\cap A}\A$,

(3) $(\tp (a,B_{1}),B_{0})$ is
a stationary pair.

\noindent
By Corollary 3.5 (v), choose finite $d\in\C$
such that $a\nda_{B_{1}}d\cup b$.
Choose $B_{2}\supseteq B_{1}\cup d$ of power
$<\k (\M )$ such that $B_{2}\subseteq\C$ and $\tp (a\cup b,\C )$ does not
split strongly over $B_{2}$.
Since $\tp (a,\C )$ and $\tp (b,\C )$ do not split strongly over
$B_{2}$ we can find by Lemmas 4.3 and 4.2 (i)
$B_{3}\supseteq B_{2}$
of power $<\k (\M )$ such that $B_{3}\subseteq\C$ and
both $(\tp (a,B_{3}),B_{2})$ and $(\tp (b,B_{3}),B_{2})$ are stationary pairs.
Then

(*) $a\da_{B_{0}}B_{3}$ and $b\da_{B_{0}\cap A}B_{3}$.

Choose $f\in Aut(B_{1})$ so that $f(B_{3})\subseteq\A$
and for all $c\in B_{3}$, $f(c)\ E^{m}_{min,B_{1}}\ c$.
Then $\tp (f(a),f(B_{3}))=\tp (a,f(B_{3}))$ and so we may assume that
$f(a)=a$.
Now $a\cup f(b)\da_{f(B_{2})}f(B_{3})$, and so we can find
$a'$ and $b'$ so that $\tp (a'\cup b',f(B_{3}))=\tp (a\cup f(b),f(B_{3}))$
and $a'\cup b'\da_{f(B_{2})}\A$. Then by (*) and
Lemma 4.2 (ii), $a'\da_{B_{0}}\A$ and so $\tp (a',\A )=\tp (a,\A )$
and we may assume that $a'=a$. Also by Lemma 4.2 (ii) and (*),
$b'\da_{B_{0}\cap A}\A$ and so $b'\da_{A}\A$. Because
$a\nda_{B_{1}}f(c)\cup b'$, by Corollary 3.5 (v),
$a\nda_{\A}b'$. Clearly this contradicts (ii). $\eop$

\th 4.7 Lemma. Let $\xi\ge\k_{r}(\M )$ be a cardinal.
Assume $D\subseteq C$,
$p\in S(C )$, $(p,D)$ is a stationary pair and
orthogonal to $\A$, $\vert C\vert <\xi$,
$\A\subseteq\B$ are strongly
$F^{\M}_{\xi}$-saturated
and
$C\da_{\A}\B$. Then $(p,D)$ is orthogonal to $\B$.

\proof For a contradiction, assume that $q\in S(\B )$ is not
orthogonal to $(p,D)$.
Choose $B\subseteq\B$ of power $<\k (\M )$
so that $q$ does not split strongly over $B$.
Choose $A\subseteq\A$ so that

(i) $\vert A\vert <\xi$,

(ii) for all $c\in C$, $\tp (c,\A\cup\B )$ does not split strongly over $A$.

\noindent
By Lemma 1.9 (v), we can find $B'\subseteq\A$
and $f\in Aut(A)$ so that $f(B)=B'$ and for all $b\in B$,
$b\ E^{m}_{min,A}\ f(b)$. By Lemma 1.12, $\tp (B',C)=\tp (B,C)$.
Let $q'=f(q)\raj B'$. Then it is easy to see that
$q'$ and $(p,C)$ are not orthogonal,
a contradiction. $\eop$

\th 4.8 Corollary. Assume $\A\subseteq\B\cap\C$ are
strongly $F^{\M}_{\k_{r}(\M )}$-saturated,
$\B\da_{\A}\C$ and $p\in S(\C )$ is orthogonal to $\A$.
Then $p$ is orthogonal to $\B$.

\proof Follows immediately from Lemma 4.3, 4.5 and 4.7. $\eop$

\chapter{5. Structure of $s$-saturated models}

We say that $\M$ is superstable if $\k (\M )=\o$.

\th 5.1 Lemma. The following are equivalent.

(i) $\k (\M )=\o$.

(ii) There are no increasing sequence
$\A_{i}$, $i<\o$,
of $a$-saturated models
and $a$ such that for all $i<\o$, $a\nda_{\A_{i}}\A_{i+1}$.

(iii) There are no increasing sequence
$\A_{i}$, $i<\o$,
of $F^{\M}_{\l (\M )}$-saturated models
and $a$ such that for all $i<\o$, $a\nda_{\A_{i}}\A_{i+1}$.

\proof Clearly (i) implies (ii)
and (ii) implies (iii). So we assume that (i)
does not hold and prove that (iii) does not hold either.
For this, choose an increasing
sequence of regular cardinals $\xi_{i}$, $i<\o$,
such that for all $i<\o$, $\M$ is $\xi_{i}$-stable.
Let $\xi =sup_{i<\o}\xi_{i}$.
By  Theorem 1.17, $\M$ is not $\xi$-stable. Let $A$
be such that $\vert A\vert\le\xi$ and $\vert S(A)\vert >\xi$.
Then choose an increasing sequence $\A_{i}$, $i<\o$,
of $F^{\M}_{\l (\M )}$-saturated models of power $\xi_{i}$
such that $A\subseteq\cup_{i<\o}\A_{i}$.
Then $\vert S(\cup_{i<\o}\A_{i})\vert >\xi$.
By Corollary 3.5 (i), it is enough to show that
there is $a$ such that
for all $i<\o$, $a\nda_{\A_{i}}\cup_{i<\o}\A_{i}$.
For a contradiction, assume not. Then for all $a$
there is $i_{a}<\o$, such that
$a\da_{\A_{i_{a}}}\cup_{i<\o}\A_{i}$. Then by Lemma 3.4,
for all $a$, $\tp (a,\cup_{i<\o}\A_{i})$ is determined by
$\tp (a,\A_{i_{a}})$. Since for all $i<\o$,
$\vert S(\A_{i})\vert\le\xi$, this implies that
$\vert S(\cup_{i<\o}\A_{i})\vert\le\xi$, a contradiction.
$\eop$

\th 5.2 Definition. We say that $\tp (a,A)$ is $F^{\M}_{\xi}$-isolated
if there is $B\subseteq A$ of power $<\xi$, such that
for all $b$, $\tp (b,B)=\tp (a,B)$ implies $\tp (b,A)=\tp (a,A)$.
We define $F^{\M}_{\xi}$-construction, $F^{\M}_{\xi}$-primary etc,
as in [Sh2]. Instead of $F^{\M}_{\l (\M )}$-saturated,
$F^{\M}_{\l (\M )}$-isolated etc, we write $s$-saturated,
$s$-isolated etc.

In slightly different context, the following theorem is proved
in [Sh1].

\th 5.3 Theorem. ([Sh1]) Assume $\xi\ge\l (\M )$.

(i) For all $A$ there is an $F^{\M}_{\xi}$-primary model
over $A$.

(ii) If $\A$ is $F^{\M}_{\xi}$-primary over $A$ then it is
$F^{\M}_{\xi}$-prime over $A$.

(iii) If $\A$ is $F^{\M}_{\xi}$-primary over $A$
and $\xi\ge\l (\M )$ is regular, then $\A$ is
$F^{\M}_{\xi}$-atomic over $A$.

(iv) If $\xi\ge\l (\M )$ is regular, then
$F^{\M}_{\xi}$-primary models over any set $A$ are unique
up to isomorphism over $A$.

As usual we write $A\do_{C}B$ if for all $a$, $a\da_{C}A$
implies $a\da_{C}B$.

\th 5.4 Lemma.

(i) Assume $\A$ is $s$-saturated and $\B$ is
$s$-primary over $\A\cup B$. Then $B\do_{\A}\B$.

(ii) Assume $\A\subseteq B\cap C$, $\A$ is
$s$-saturated and $B\da_{\A}C$. If
$(B,\{ b_{i}\vert\ i<\g\} ,(B_{i}\vert\ i<\g ))$ is an
$s$-construction over $B$, then
$(B\cup C,\{ b_{i}\vert\ i<\g\} ,(B_{i}\vert\ i<\g ))$
is an $s$-construction over $B\cup C$.

(iii) Assume $\xi\ge\l (\M )$, $\A$ is $F^{\M}_{\xi}$-saturated and
$\B$ is $F^{\M}_{\xi}$-primary over $\A\cup B$. Then
$B\do_{\A}\B$.

\proof (i): Assume not. Then we can find
$s$-saturated $\A$, $B$,
$b$ and $a$ so that $\tp (b,\A\cup B)$ is $s$-isolated,
$a\da_{\A}B$ and $a\nda_{\A}b$
(if $(\A\cup B,\{ b_{i}\vert\ i<\g\} ,(B_{i}\vert\ i<\g ))$
is an $s$-construction of $\B$, then let $b=b_{i}$,
where $i$ is the least ordinal such
that $a\nda_{\A}B\cup\bigcup_{j\le i}b_{j}$
and rename $B\cup\bigcup_{j<i}b_{j}$ as $B$; $i$ exists by
Corollary 3.5 (v)).
Without loss of generality
we may assume that $\vert B\vert <\l (\M )$. Choose
$A\subseteq\A$ so that

(i) $\tp (b,A\cup B)$ $s$-isolates $\tp (b,\A\cup B)$,

(ii) for all $c\in B$, $\tp (c,\A\cup a)$ does not split strongly over
some $A'\subseteq A$ of power $<\k (\M )$,

(iii) $\tp (b,\A )$ does not split strongly over
some $A'\subseteq A$ of power $<\k (\M )$,

(iv) $\vert A\vert <\l (\M )$.

\noindent
This is possible since $\k_{r}(\M )\le\l (\M)$:
Let $\d =\vert B\vert +1 <\l (\M )$.
Clearly we can choose $A$ so that it
of the form $A'\cup A''$ where $A'$ is of power $<\l (\M )$ and $A''$
is a union of $\d$ many sets of power $<\k_{r}(\M )\le\l (\M )$.
If $\l (\M )$ is regular, then clearly $\vert A\vert <\l (\M )$.
Otherwise $\k_{r}(\M )<\l (\M )$ in which case
$\vert A\vert\le\vert A'\vert +max(\d ,\k_{r}(\M ))<\l (\M )$.

By Lemma 1.9 (iii),
the proof of Lemma 1.13 and (iii) above, there are $c,c',a'\in\A$ such that
$c\cup a\ E^{m}_{min,A}\ c'\cup a'$ and
$\tp (b\cup c\cup a,A)\ne \tp (b\cup c'\cup a',A)$. By (ii),
$\tp (B\cup c\cup a,A)=\tp (B\cup c'\cup a',A)$. So
there is $f\in Aut(A\cup B)$ such that $f(c)=c'$ and
$f(a)=a'$. Then $f(b)$ contradicts (i) above.

(ii): As (i) above.

(iii): By (i) we may assume that $\xi >\l (\M )$.
For a contradiction, assume that the claim does not hold.
As in (i), we can find
$s$-saturated $\A$, $B$,
$b$ and $a$ so that $\tp (b,\A\cup B)$ is $F^{\M}_{\xi}$-isolated,
$a\da_{\A}B$, $a\nda_{\A}b$ and
$\vert B\vert <\xi$.
Let $A\subseteq\A$ be such that $\tp (a,A\cup B)$
$F^{\M}_{\xi}$-isolates $\tp (b,\A\cup B)$.
Choose $s$-saturated $\C\subseteq\A$
so that $\vert\C\vert =\l (\M )$ and
$a\da_{\C}\A\cup\B$. For $i<\xi$, choose $a_{i}\in\A$ such that
$(a_{i})_{i<\xi}$ is $\C$-independent and for all $i<\xi$,
$\tp (a_{i},\C )=\tp (a,\C )$. As in (i), it is enough to
show that there is $i<\xi$ such that $a_{i}\da_{\C}A\cup B$.
For this we choose maximal sequence of models $\A_{j}$
and sets $I_{j}\subseteq\xi$,
$j\le j^{*}$, such that

(a) $\A_{0}=\C$ and $I_{0}=\empty$,

(b) $I_{j+1}-I_{j}$ is finite, $\A_{j+1}$ is
$s$-primary over $\A_{j}\cup (I_{j+1}-I_{j})$ and
for some $c\in A\cup B$, $c\nda_{\A_{j}}I_{j+1}-I_{j}$,

(c) if $j$ is limit, then $I_{j}=\cup_{k<j}I_{k}$ and
$\A_{j}$ is $s$-primary over $\cup_{k<j}\A_{k}$.

\noindent
Since $\k_{r}(\M )\le\vert A\cup B\vert <\xi$,
$I_{j^{*}}\ne\xi$. Let $i\in\xi -I_{j^{*}}$.
By (i) and (ii), it is easy to see that for all $j\le j^{*}$,
$\A_{j}$ is $s$-primary over $\A\cup I_{j}$. Then
by (i), $a_{i}\da_{\C}\A_{j^{*}}$ and because
the sequence was maximal, $A\cup B\da_{\A_{j^{*}}}a_{i}$.
So $a_{i}\da_{\C}A\cup B$ as wanted.
$\eop$

\th 5.5 Corollary.

(i) Assume $A\subseteq\A$ and $\A$ is
$s$-saturated. If $p\in S(\A )$ is orthogonal to $A$,
then for all $C\supseteq\A$, $a$ and $b$ the following holds:
if $a\da_{\A}C$, $\tp (a,\A )=p$ and $b\da_{A}C$, then
$a\da_{\A}C\cup b$.

(ii) Assume $\M$ is superstable and $\g$ is a limit ordinal.
Let $\A_{i}$, $i<\g$, be an increasing sequence of
$s$-saturated models and $\A$ be $s$-primary over $\cup_{i<\g}\A_{i}$.
If $a\not\in\A$ then there is $i<\g$ such that
$\tp (a,\A )$ is not orthogonal to $\A_{i}$.

(iii) Assume $\A$ is $s$-saturated and $p\in S(\A )$ is orthogonal
to $A\subseteq\A$.
If $a_{i}$, $i<\o$, are such that for all $i<\o$,
$\tp (a_{i},\A )=p$ and $a_{i}\da_{\A}\cup_{j<i}a_{j}$, then
for all $n<\o$, $\tp (\cup_{i<n}a_{i},\A )$
is orthogonal to $A$.

\proof (i): Assume not. Let $\C$ be $s$-primary over $\A\cup C$.
Then by Lemma 5.4 (i) and Corollary 3.5 (iv),
$a\da_{\A}\C$, $b\da_{A}\C$ and $a\nda_{\C}b$,
a contradiction.

(ii): Clearly we may assume that if $i<j$ then $\A_{i}\ne\A_{j}$.
Since $\k (\M )=\o$, there is $i<\g$ such that
$a\da_{\A_{i}}\cup_{j<\g}\A_{j}$. By (i),
$a\da_{\A_{i}}\A$. By Lemma 3.2 (v), this is more that required.

(iii): Assume not. Then by Lemma 4.6, there is $b$ such that
$b\da_{A}\B$ and $\cup_{i<n}a_{i}\nda_{\A}b$. Let $m\le n$
be the least such that $\cup_{i<m}a_{i}\nda_{\A}b$. By Lemma 3.8 (i),
$a_{m-1}\nda_{\A\cup\bigcup_{i<m-1}a_{i}}b$. Clearly this contradicts
(i). $\eop$

Let $P$ be a tree without branches of length
$>\o$. Then
by $t^{-}$ we mean the
immediate predecessor of $t$ if $t\in P$ is not the root.
For all $t\in P$, by $t^{1}_{>}$ we mean the set of
immediate successors of $t$.

\th 5.6 Definition. ([Sh2]) We say that $(P,f,g)=((P,\prec ),f,g)$
is an $s$-free tree of $s$-saturated
$\A$
if the following holds:

(i) $(P,\prec )$ is a tree without branches of length $>\o$,
$f:(P-\{ r\} )\rightarrow\A$ and $g:P\rightarrow P(\A )$,
where $r\in P$ is the root of $P$ and $P(\A )$ is the power set
of $\A$ - in order to simplify the notation we write $a_{t}$ for
$f(t)$ and $\A_{t}$ for $g(t)$,

(ii) $\A_{r}$ is $s$-primary model (over $\empty$),

(iii)  if $t$ is not the root and $u^{-}=t$ then
$\tp (a_{u},\A_{t})$ is orthogonal to $\A_{t^{-}}$,

(iv) if $t=u^{-}$ then $\A_{u}$ is $s$-primary over $\A_{t}\cup a_{u}$,

(v) Assume $T,V\subseteq P$ and $u\in P$ are such that

(a) for
all $t\in T$, $t$ is comparable with $u$,

(b) $T$ is downwards closed.

(c) if $v\in V$ then for all $t$ such that $v\succeq t\succ u$, $t\not\in T$.

\noindent
Then
$$\bigcup_{t\in T}\A_{t}\da_{\A_{u}}\bigcup_{v\in V}\A_{v}.$$

\th 5.7 Definition. We say that $(P,f,g)$ is
an $s$-decomposition of $\A$ if it is a maximal
$s$-free tree of
$\A$.

Notice that 'the finite character of dependence' implies, that
unions of increasing sequences of $s$-free trees of $\A$ are
$s$-free trees of $\A$.
So
for all $s$-saturated $\A$ there is
an $s$-decomposition of $\A$.

We say that $\A$ is $s$-primary over an $s$-free tree
$(P,f,g)$
if $\A$ is $s$-primary over $\bigcup\{\A_{t}\vert\ t\in P\}$.

\th 5.8 Definition. Assume that $(P,f,g)$ is an $s$-decomposition of
$\A$, $\A$ is $s$-saturated. Let $P=\{ t_{i}\vert\ i<\a\}$ be
an enumeration of $P$ such that if $t_{i}\prec t_{j}$ then $i<j$. Then
we say that $(\A_{i})_{i\le\a}$ is a generating sequence if
the following holds

(i) for all $i\le\a$, $\A_{i}\subseteq\A$,

(ii) $\A_{0}=\empty$,

(iii) $\A_{i+1}$ is $s$-primary over $\A_{i}\cup\A_{t_{i}}$,

(iv) if $0<i\le\a$ is limit then $\A_{i}$ is $s$-primary over
$\bigcup_{j<i}\A_{j}$.

\th 5.9 Lemma. Assume that $(P,f,g)$ is an $s$-free tree of
$\A$, $\A$ is $s$-saturated and
$(\A_{i})_{i\le\a}$ is a generating sequence. Then for all
$0<i<\a$, $\A_{t_{i}}\da_{\A_{t_{i}^{-}}}\A_{i}$.

\proof By Lemma 5.4 (i), it is enough to prove that for all $i<\a$,
$\A_{i}$ is $s$-primary over $\cup_{j<i}\A_{t_{j}}$. We prove this
by induction on $i$. In fact we need to prove slightly more
to keep the induction going: We show that
$\A_{i}$ is not only $s$-constructible over $\cup_{j<i}\A_{t_{j}}$
but that the natural construction works. Then the limit cases are
trivial and the successor cases follow from Lemma 5.4 (ii).
$\eop$

\th 5.10 Definition. Assume $\A$ is $s$-saturated.
We say that $\tp (a,\A )$ is a c-type if for all
$s$-saturated $\C$ and $\B$ the following holds:
If $\C\subseteq\A$ is such that $\tp (a,\A )$ is not
orthogonal to $\C$ and $\A\cup a\subseteq\B$, then
there is $b\in\B -\A$ such that $b\da_{\C}\A$.

Notice that the notion of c-type is a generalization of regular type.

\th 5.11 Lemma. Assume $\M$ is superstable.
Let $\A\subseteq\B$ be $s$-saturated
and $\A\ne\B$. Then there is a singleton $a\in\B -\A$
such that $\tp (a,\A )$ is a $c$-type.

\proof Since $\k (\M )=\o$, by Lemma 1.1 it is easy to see
that there is a singleton $a\in\B -\A$ and finite $A\subseteq\A$
such that the following holds: for all $b\in\B -\A$ and $B\subseteq\A$,
if there is an automorphism $f$ of $\M$ such that $f(a)=b$ and
$f(A)=B$, then $\tp (b,\A )$ does not split strongly over $B$
(and so $b\da_{B}\A$).
We show that $a$ is as wanted. Let $s$-saturated
$\C\subseteq\A$ be such that $\tp (a,\A )$
is not orthogonal to $\C$.
Since $\B$ can now be any $s$-saturated model such that
$\A\cup a\subseteq\B$, it is enough to show that there is
$b\in\B -\A$ such that $b\da_{\C}\A$.

By Lemma 4.6, find $d$ such that $d\da_{\C}\A$ and $a\nda_{\A}d$.
Let $\D$ be $s$-primary over $\C\cup d$. Then $\D\da_{\C}\A$
and $a\nda_{\A}\D$. For all $i<\o$, choose $\A_{i}$ and $a_{i}$
so that $\tp (a_{i}\cup\A_{i},\D )=\tp (a\cup\A ,\D )$ and
$a_{i}\cup\A_{i}\da_{\D}a\cup\A\cup\bigcup_{j<i}(a_{j}\cup\A_{j})$.

{\bf Claim.} $\{ a\cup\A\}\cup\{ a_{i}\cup\A_{i}\vert\ i<\o\}$
is indiscernible
over $\C$ and $a\cup\A\nda_{\C}\cup_{i<\o}(a_{i}\cup\A_{i})$.

\proof The first of the claims follow immediately from Corollary 3.5 (ii).
For a contradiction, assume that the second claim is not true.
For all $i<\o$, we define $\B_{i}$ as follows:
We let $\B_{0}$ be $s$-primary over $\A\cup a$ and
$\B_{i+1}$ be $s$-primary over $\B_{i}\cup\A_{i}\cup a_{i}$.
By Lemma 5.1, there is $i<\o$ such that
$d\da_{\B_{i}}\A_{i}\cup a_{i}$. Since
$\{ a\cup\A\}\cup\{ a_{i}\cup\A_{i}\vert\ i<\o\}$
is indiscernible
over $\C$ and $a\cup\A\da_{\C}\cup_{i<\o}(a_{i}\cup\A_{i})$,
$\A_{i}\cup a_{i}\da_{\C}
a\cup\A\cup\bigcup_{j<i}(a_{j}\cup\A_{j})$. By Lemma 5.4 (ii),
$\A_{i}\cup a_{i}\da_{\C}\B_{i}$. But then
$\A_{i}\cup a_{i}\da_{\C}d$, a contradiction. $\eop$ Claim.

By Claim and Corollary 3.5 (v), let $n<\o$ be the least such that
$a\cup\A\nda_{\C}\cup_{i<n}(a_{i}\cup\A_{i})$. Let $\A^{*}$
be $s$-primary over $\A\cup\A_{0}\cup\bigcup_{0<i<n}(\A_{i}\cup a_{i})$.
It is easy to see that
$\A_{n}\da_{\C}\A\cup\bigcup_{0<i<n}(\A_{i}\cup a_{i})$. By Claim,
$\A_{0}\da_{\C}\A\cup\bigcup_{0<i<n}(\A_{i}\cup a_{i})$ and
so by Lemma 3.8 (iv) and the choice of $n$,
$\A\cup a\da_{\C}\A_{0}\cup\bigcup_{0<i<n}(\A_{i}\cup a_{i})$
and so by Lemma 3.6 and 3.2 (i),
$a\da_{\A}\A_{0}\cup\bigcup_{0<i<n}(\A_{i}\cup a_{i})$.
By Lemma 5.4 (i), $a\da_{\A}\A^{*}$.
Similarly we we see that
$a_{0}\da_{\A_{0}}\A^{*}$.
Then also $a\nda_{\A^{*}}a_{1}$.

By the choice of $\A_{0}$ and $a_{0}$ there is $f\in Aut(\C )$
such that $f(\A )=\A_{0}$ and $f(a)=a_{0}$.
Let $A_{0}=f(A)$. By Corollary 3.5 (v) there is finite
$C\subseteq\A^{*}$ such that
$a\nda_{A}A_{0}\cup C\cup a_{0}$. Choose $B\subseteq\C$ such that
$\tp (A\cup a,\C)$ does not split strongly over $B$.
Then there is $g\in Saut(B)$ such that $g(A_{0})\subseteq\C$.
Since $a\cup\A\da_{\C}\A^{*}$ and
every $h\in Aut(\A^{*})$ belongs to $Saut(B)$, we may assume that
$$(*)\ \ \ \ \ a\cup\A\da_{\C}g(C)\cup A_{0}\cup C.$$
Then
$\tp (g(A_{0}\cup C),A\cup a)=\tp (A_{0}\cup C,A\cup a)$.
Choose $h\in Saut(A\cup g(A_{0}))$ such that
$h(g(C))\subseteq\A$. By (*),
$\tp (a,\A\cup g(C))$ does not split strongly over $A$ and so
it does not split strongly
over $A\cup g(A_{0})$. Then
$\tp (g(A_{0})\cup h(g(C)),A\cup a)=\tp (A_{0}\cup C,A\cup a)$.
Choose $b\in\B$ such that
$\tp (g(A_{0})\cup h(g(C))\cup b,A\cup a)=\tp (A_{0}\cup C\cup a_{0},A\cup a)$.
Then by Corollary 3.5 (v) and the choice of $C$,
$a\nda_{\A}b$ and so by Lemma 3.2 (iii), $b\in\B -\A$
($b$ is a singleton). By the choice of $A$,
$\tp (b,\A )$ does not split strongly over $g(A_{0})$.
By Lemma 3.2 (iii), $b\da_{\C}\A$. $\eop$

\th 5.12 Definition.

(i) We say that $\M$ has $s$-SP (structure property) if every
$s$-saturated $\A$ is $s$-primary over any $s$-decomposition of $\A$.

(ii) Let $\k\ge\l (\M )$.
We say that $\M$ has $\k$-dop if there are $F^{\M}_{\k}$-saturated
$\A_{i}$, $i<4$, and $a\not\in\A_{3}$ such that

(a) $\A_{0}\subseteq\A_{1}\cap\A_{2}$, 

(b) $\A_{1}\da_{\A_{0}}\A_{2}$,

(c) $\A_{3}$ is $F^{\M}_{\k}$-primary over $\A_{1}\cup\A_{2}$,

(d) $\tp (a,\A_{3})$ is orthogonal to $\A_{1}$
and to $\A_{2}$.

\noindent
We say that $\M$ has $\k$-ndop if it does not have $\k$-dop.

\th 5.13 Theorem. Assume $\M$ is superstable and has $\l (\M )$-ndop.
Then $\M$ has $s$-SP.

\proof
Let $\A$ be $s$-saturated and $(P,f,g)$ an
$s$-decomposition of $\A$.
Let $(\A_{i})_{i\le\a}$ be a
generating sequence and $P=\{ t_{i}\vert\ i<\a\}$ be the
enumeration of $P$ from the definition of a generating sequence.

{\bf Claim}: $\A_{\a}=\A$.

\proof Assume not. For all $a\in\A -\A_{\a}$ let $i_{a}$ be the least
ordinal such that $\tp (a,\A_{\a})$ is not orthogonal to $\A_{i_{a}}$. Let
$a\in\A -\A_{\a}$ be any sequence such that

(i) for some
$l\le\a$ either $\tp (a,\A_{l} )$ is a $c$-type
and $a\da_{\A_{l}}\A_{\a}$ or $\tp (a,\A_{t_{l}} )$ is a $c$-type
and $a\da_{\A_{t_{l}}}\A_{\a}$

\noindent
and

(ii) among these $a$, $i=i_{a}$ is the least.

\noindent
By Lemma 5.11 there is at least one such $a$.

There are two cases:

Case 1: For some $l<\a$
$\tp (a,\A_{t_{l}} )$ is a c-type
and $a\da_{\A_{t_{l}}}\A_{\a}$.
Let $t^{*}\le t_{l}$ be the least $t$ such that
$\tp (a,\A_{t_{l}} )$ is not orthogonal to $\A_{t}$. Since
$\tp (a,\A_{t_{l}} )$ is a $c$-type choose $b$ so that

(1) $b\da_{\A_{t^{*}}}\A_{t_{l}}$

\noindent
and

(2) $b\in\A_{t_{l}}[a] -\A_{t_{l}}$, where $\A_{t_{l}}[a]\subseteq\A$
is $s$-primary over $\A_{t_{l}}\cup a$.

\noindent
Then if $(t^{*})^{-}$ exists, by (2) and Lemmas 4.6 and 5.4 (i),
$\tp (b,\A_{t_{l}})$ is orthogonal to $\A_{(t^{*})^{-}}$
and so by (1) and Lemma 4.6 it is
easy to see that $\tp (b,\A_{t^{*}})$ is orthogonal to $\A_{(t^{*})^{-}}$.

By (1), (2) and Lemma 5.4 (i),
$b\da_{\A_{t^{*}}}\A_{\a}$.

We define $((P',\prec'),f',g')$ as follows:

(i) $P'=P\cup\{ t\}$, $t$ a new node,

(ii) for all $u\in P$, $u\prec't$ iff $u\preceq t^{*}$

(iii) $f'\raj P=f$ and $f'(t)=b$,

(iv) $g'\raj P=g$ and $g'(t)\subseteq\A$ is $s$-primary
over $A_{t^{*}}\cup b$.

{\bf Subclaim.} $((P',\prec'),f',g')$ is an $s$-free tree of $\A$.

\proof (i), (ii), (iii) and (iv) in the Definition 5.6 are clear.
So we prove (v):

Let $T\subseteq P'$, $u\in P'$ and $V\subseteq P'$ be as
in Definition 5.6 (v).
There are four cases:

Case a: $t\in T-V$.
Let $T'=T-\{ t\}$ and $\A_{T'}\subseteq\A_{\a}$ be $s$-primary
over $\cup\{\A_{d}\vert\ d\in T'\}$.
By the choice of $b$ and Lemma 5.4 (i),
$$A_{t}\da_{\A_{t^{*}}}\A_{T'}\cup\bigcup_{v\in V}\A_{v}.$$
By Lemmas 3.2 (i) and 3.6,
$$\cup_{v\in V}\A_{v}\da_{\A_{T'}}\A_{t}.$$
By Corollary 3.5 (iv), the assumption that
$(P,f,g)$ is $s$-free tree of $\A$ and Lemma 5.4 (i),
$$\cup_{v\in V}\A_{v}\da_{\A_{u}}\A_{T'}\cup\A_{t}.$$
By Lemma 3.6,
$$\cup_{d\in T}\A_{d}\da_{\A_{u}}\cup_{v\in V}\A_{v}.$$

Case b: $t\in V-T$: Exactly as the Case a.

Case c: $t\in V\cap T$: Because $t\in T-P$, $u\le t$. Since
$t\in V$, $u=t$. Then because $u\not\in P$,
$\cup_{d\in T}\A_{d}=\A_{u}$, and the
claim follows from Lemma 3.2 (iv).

Case d: $t\not\in T\cup V$: Immediate by the assumption that
$(P,f,g)$ is an $s$-free tree of $\A$.

$\eop$ Subclaim.

Subclaim contradicts the maximality of $P$. So Case 1 is impossible
and we are in the Case 2:

Case 2: $l\le\a$ is such that
$\tp (a,\A_{l} )$ is a $c$-type
and $a\da_{\A_{l}}\A_{\a}$.
Let $\B\subseteq\A$ be $s$-primary over $\A_{\a}\cup a$.
Clearly $i(=i_{a})\le l$ and so
let $b'$ be the element given by $\tp (a,\A_{l} )$
being a $c$-type: $b'\da_{\A_{i}}\A_{l}$ and
$b'\in\A_{l}[a]-\A_{l}$, where $\A_{l}[a]\subseteq\B$ is $s$-primary over
$\A_{l}\cup a$.
By Lemma 5.11 we may choose $b$ so that
$\tp (b,\A_{i})$ is a $c$-type and $b\in\A_{i}[b']-\A_{i}$, where
$\A_{i}[b']\subseteq\B$ is $s$-primary over
$\A_{i}\cup b'$.
Then
$b\da_{\A_{i}}\A_{\a}$, $b\not\in\A_{i}$ and $i_{b}\le i(=i_{a})$.

1. $i$ is not a limit $>0$. This is because otherwise by
Lemma 5.5 (ii), $\tp (b,\A_{i})$ is not orthogonal to
$\A_{j}$ for some $j<i$.
Then $\tp (b,\A_{\a})$ is not orthogonal to
$\A_{j}$ i.e. $i_{b}<i_{a}$.
This contradicts the choice of $a$.

2. $i$ is not a successor $>1$. Assume it is, $i=j+1$. Then
$\A_{i}$ is $s$-primary over $\A_{j}\cup\A_{t_{j}}$ and
by Lemma 5.9,
$\A_{j}\da_{\A_{t_{j}^{-}}}\A_{t_{j}}$.
(Notice that since Case 1 is not possible, $A_{j+1}\ne\A_{t_{j}}$.)
By the choice of $a$
$\tp (b,\A_{i})$ is orthogonal to $\A_{j}$.
So by $\l (\M )$-ndop $\tp (b,\A_{i})$ is not orthogonal to $\A_{t_{j}}$.
Then as in Case 1 we get a contradiction with the maximality
of $(P,f,g)$. Alternatively, we can find $c$ such that it
satisfies the assumptions of Case 1, which is
a contradiction.

3. $i$ is not $0$ or $1$. Immediate, since Case 1 is not possible.

Clearly 1 and 2 above contradict 3.
So also Case 2 imply a contradiction.

$\eop$ Claim.

Let $\C\subseteq\A$ be $F^{\M}_{\k}$-primary
over $\bigcup\{ \A_{t}\vert\ t\in P\}$.
We want to show that $\C =\A$. For this we choose
a generating sequence $(\A_{i})_{i\le\a}$,
so that
$\A_{i}\subseteq\B$ for all $i\le\a$. By the claim above
$\A_{\a}=\A$ and so $\C =\A$.
$\eop$

\chapter{6. On nonstructure}

\th 6.1 Definition. We say that $\M$ has $\k$-sdop if the
following holds: there are $F^{\M}_{\k}$-saturated
$\A_{i}$, $i<4$, and
$I=\{ a_{i}\vert\ i<\l (\M )\}$, $a_{i}\in\A_{3}$, such that

(a) $\A_{0}\subseteq\A_{1}\cap\A_{2}$, $\A_{3}$
is $F^{\M}_{\k}$-primary over $\A_{1}\cup\A_{2}$,

(b) $\A_{1}\da_{\A_{0}}\A_{2}$,

(c) $I$ is an indiscernible sequence over $\A_{1}\cup\A_{2}$
and if $i<j<\l (\M )$ then $a_{i}\ne a_{j}$.

As in [Hy],
we can prove non-structure theorems from $\k$-sdop.
(In [Sh2], this was the formulation of dop, which was used to get
non-structure.)

In this chapter we show that dop and sdop are essentially
equivalent i.e. $\l (\M )^{+}$-sdop implies $\l (\M )^{+}$-dop
and $\l (\M )$-dop implies $\l_{r}(\M )^{+}$-sdop, where
$\l_{r}(\M )$ is the least regular cardinal $\ge\l (\M )$.

\th 6.2 Lemma. Assume $\M$ is $\xi$-stable and $\k =\xi^{+}$.
If $\M$ has $\k$-sdop then it has $\k$-dop.

\proof Let $I$ and $\A_{i}$, $i<4$, be as in the definition of
$\k$-sdop. We need to show
that there is $\M$-consistent type $p$ over $\A_{3}$ such that (d)
in Definition 5.12 (ii) is satisfied.
We show that $Av(I,\A_{3})$
is the required type.

By Lemma 2.4 (iii),
let $a$ be such that $\tp (a,\A_{3})=Av(I,\A_{3})$.
For a contradiction, by Lemma 4.6,
let $b$ be such that

(i) $a\nda_{\A_{3}}b$,

(ii) $b\da_{\A_{1}}\A_{3}$.

\noindent
Let $\C_{i}\subseteq\A_{i}$, $i<4$ be $F^{\M}_{\xi}$-saturated
models of cardinality
$\xi$ such that

(1) $\C_{i}\subseteq\A_{i}$, $\C_{1}\cap\C_{2}=\C_{0}$,
$\C_{3}\cap\A_{1}=\C_{1}$, $\C_{3}\cap\A_{2}=\C_{2}$ and
$I\subseteq\C_{3}$,

(2) $a\cup b\da_{\C_{3}}\A_{3}$ and $a\nda_{\C_{3}}b$,

(3) $a\cup b\cup\C_{3}\da_{\C_{1}}\A_{1}$ and
$a\cup b\cup\C_{3}\da_{\C_{2}}\A_{2}$,

(4) for all $c\in\C_{3}$ there is $D\subseteq\C_{1}\cup\C_{2}$
of power $\xi$,
such that $\tp (c,D)$ $F^{\M}_{\k}$-isolated
$\tp (c,\A_{1}\cup\A_{2})$.

We can see the existence of the sets as in the proof of Theorem 3.14
(the only non-trivial part being to guarantee that the
models are $F^{\M}_{\xi}$-saturated).

Let $a^{*}\in\A_{3}$ be such that it realizes $Av(I,\C_{3})$.

{\bf Claim.} $\tp (a^{*},\C_{3})$
$F^{\M}_{\k}$-isolates $\tp (a^{*},\C_{3}\cup\A_{1}\cup\A_{2})$.

\proof Assume not. Then there is $d\in\C_{3}$ such that
$\tp (a^{*}\cup d,\C_{1}\cup\C_{2})$ does not
$F^{\M}_{\k}$-isolate $\tp (a^{*}\cup d,\A_{1}\cup\A_{2})$.

{\bf Subclaim.} There is $a'\in I$ such that
$\tp (a'\cup d,\C_{1}\cup\C_{2})=\tp (a^{*}\cup d,\C_{1}\cup\C_{2})$.

\proof By Lemma 1.2 (v), there is $i<\l (\M )$ such that
$\tp (d,\C_{1}\cup\C_{2}\cup I)$ does not split over
$\C_{1}\cup\C_{2}\cup\{ a_{j}\vert\ j<i\}$. Since
$I$ is indiscernible over $C_{1}\cup\C_{2}$, $a'=a_{i}$
is as wanted. $\eop$ Subclaim.

Clearly Subclaim contradicts (4) above. $\eop$ Claim.

Choose $b^{*}\in\A_{1}$ so that $\tp (b^{*},\C_{1})=\tp (b,\C_{1})$.
By (3), $\tp (b^{*},\C_{3})=\tp (b,\C_{3})$. By Claim,
$a^{*}\da_{\C_{3}}b^{*}$. Let $f$ be an automorphism such that
$f(b)=b^{*}$ and $f\raj\C_{3}=id_{\C_{3}}$. Then $f(a)$ contradicts
Claim. $\eop$

\th 6.3 Theorem. Let $\l\ge\l_{r}(\M )$ be such that $\M$ is
$\l$-stable. Then
$\l_{r}(\M )$-dop implies
$\l^{+}$-sdop.

\proof Let $\A_{i}$, $i<4$, and $p\in S(\A_{3})$ be as in the
the definition of $\l_{r}(\M )$-dop. By Lemma 4.5,
as in the proof of Lemma 3.14, we find these so
that $\vert\A_{3}\vert\le\l_{r}(\M )$.
Let $\B_{0}\supseteq\A_{0}$ be $F^{\M}_{\l^{+}}$-saturated
such that $\B_{0}\da_{\A_{0}}\A_{3}$. Let $\B_{i}$, $i\in\{ 1,2\}$
be $s$-primary over $\A_{i}\cup\B_{0}$.
Let $\B_{3}$ be
$s$-primary over $\B_{1}\cup\B_{2}\cup\A_{3}$.
Let $\C_{i}$, $i\in\{ 1,2\}$, be $F^{\M}_{\l^{+}}$-primary
over $\B_{i}$ such that $\C_{1}\da_{\B_{1}}\B_{3}$ and
$\C_{2}\da_{\B_{2}}\B_{3}\cup\C_{1}$.

Let $q\in S(\A_{3})$ be any type such that it is orthogonal
to $\A_{1}$ and $\A_{2}$. Our first goal is to show that
there is only one $q^{*}\in S(\C_{1}\cup\C_{2}\cup\B_{3})$
which extends $q$.

{\bf Claim 1.} $\B_{1}\cup\A_{2}$ is $F^{\M}_{\l_{r}(\M )}$-constructible over
$\A_{1}\cup\B_{0}\cup\A_{2}$ and for all $b\in\B_{1}$, there
is $B\subseteq\A_{1}\cup\B_{0}$ of power $<\l_{r}(\M )$ such that
$\tp (b,B)$ $F^{\M}_{\l_{r}(\M )}$-isolates $\tp (b,\B_{0}\cup\A_{1}\cup\A_{2})$.

\proof Follows immediately from the proof of Lemma 5.4 (ii)
and Theorem 5.3 (iii). $\eop$ Claim 1.

{\bf Claim 2} $\B_{1}\cup\B_{2}$ is $F^{\M}_{\l_{r}(\M )}$-constructible over
$\B_{1}\cup\B_{0}\cup\A_{2}$ and for all $b\in\B_{2}$, there
is $B\subseteq\A_{2}\cup\B_{0}$ of power $<\l_{r}(\M )$ such that
$\tp (b,B)$ $F^{\M}_{\l_{r}(\M )}$-isolates $\tp (b,\B_{1}\cup\A_{2})$.

\proof As Claim 1. $\eop$ Claim 2.

{\bf Claim 3} $\B_{1}\cup\B_{2}\cup\A_{3}$ is
$F^{\M}_{\l_{r}(\M )}$-constructible over $\A_{3}\cup\B_{0}$.

\proof By Claims 1 and 2, $\B_{1}\cup\B_{2}$ is
$F^{\M}_{\l_{r}(\M )}$-constructible
over $\A_{1}\cup\A_{2}\cup\B_{0}$. So it is enough to show
that for all $a\in\A_{3}$, $\tp (a,\A_{1}\cup\A_{2})$
$F^{\M}_{\l_{r}(\M )^{+}}$-isolates $\tp (a,\B_{1}\cup\B_{2})$.

Assume not. Choose $b_{1}\in\B_{1}$ and $b_{2}\in\B_{2}$
so that $\tp (a,\A_{1}\cup\A_{2})$ does not
$F^{\M}_{\l_{r}(\M )^{+}}$-isolate $\tp (a,\A_{1}\cup\A_{2}\cup b_{1}\cup b_{2})$.
Choose $A_{1}\subseteq\A_{1}$, $A_{2}\subseteq\A_{2}$ and
$B_{0}\subseteq\B_{0}$
of power $<\l_{r}(\M )$ such that

(i) $\tp (a,A_{1}\cup A_{2})$
$F^{\M}_{\l_{r}(\M )}$-isolates $\tp (a,\A_{1}\cup\A_{2})$,

(ii) $\tp (b_{1},A_{1}\cup B_{0})$ $F^{\M}_{\l_{r}(\M )}$-isolates
$\tp (a,\A_{1}\cup\A_{2}\cup\B_{0})$ and
$\tp (b_{2},A_{2}\cup B_{0})$ $F^{\M}_{\l_{r}(\M )}$-isolates
$\tp (b_{2},\B_{1}\cup\A_{2}\cup\B_{0})$,

(iii) $A_{1}\cap\A_{0}=A_{2}\cap\A_{0}=A_{0}$ and for all
$c\in A_{1}\cup A_{2}$, $\tp (c,\B_{0})$ does not split strongly
over $A_{0}$.

\noindent
By Lemma 1.9 (v), choose $f\in Aut(A_{0})$ so that
$f(B_{0})\subseteq\A_{0}$ and for all $c\in B_{0}$,
$f(c)\ E^{m}_{min,A_{0}}\ c$.
Let $B'_{0}=f(B_{0})$.
Then by (iii), $\tp (B'_{0},A_{1}\cup A_{2})=\tp (B_{0},A_{1}\cup A_{2})$.
Choose $b'_{i}\in\A_{i}$ so that
$\tp (b'_{i}\cup B'_{0},A_{i})=\tp (b_{i}\cup B_{0},A_{i})$, $i\in\{ 1,2\}$.
By (ii) $\tp (b'_{1}\cup b'_{2}\cup B'_{0},A_{i})=
\tp (b_{1}\cup b_{2}\cup B_{0},A_{i})$. Clearly this contradicts
(i). $\eop$ Claim 3.

{\bf Claim 4.} $\B_{3}$ is
$F^{\M}_{\l_{r}(\M )}$-primary over $\A_{3}\cup\B_{0}$.

\proof Immediate by Claim 3 and the choice of $\B_{3}$. $\eop$ Claim 4.

By Claim 4 and Lemma 5.4, there is exactly one $q'\in S(\B_{3})$
such that $q\subseteq q'$. By Corollary 4.8,
$q'$ is orthogonal to $\B_{1}$ and $\B_{2}$.
So if $a$ realizes $q'$, then $a\da_{\B_{3}}\C_{1}$.
Then by Corollary 5.5 (i), there is exactly one
$q^{*}\in S(\C_{1}\cup\C_{2}\cup\B_{3})$,
which extends $q$.

Now choose $a_{i}$, $i<\l_{r}(\M )$, so that for all
$i$, $\tp (a_{i},\A_{3})=p$ and $a_{i}\da_{\A_{3}}\cup_{j<i}a_{j}$.
Then $I=\{ a_{i}\vert\ i<\l_{r}(\M )\}$ is
indiscernible over $\A_{3}$ and by Corollary 5.5 (iii), for all $n<\o$,
$\tp (a_{0}\cup ...\cup a_{n},\A_{3})$ is orthogonal to
$\A_{1}$ and $\A_{2}$. So, by what we showed above,
$I$ is indiscernible over $\C_{1}\cup\C_{2}\cup\B_{3}$ and
for all $i<\l_{r}(\M )$,
$\tp (a_{i},\A_{3}\cup\bigcup_{j<i}a_{j})$ $F^{\M}_{\l^{+}}$-isolates
$\tp (a_{i},\C_{1}\cup\C_{2}\cup\B_{3}\cup\bigcup_{j<i}a_{j})$.
So there is an $F^{\M}_{\l^{+}}$-primary model $\C_{3}$ over
$\C_{1}\cup\C_{2}\cup\B_{3}$ such that $I\subseteq\C_{3}$.

So to get $\l^{+}$-sdop, it is enough to show that
$\C_{3}$ is $F^{\M}_{\l^{+}}$-primary over $\C_{1}\cup\C_{2}$.
By Claim 3 and the choice of $\B_{3}$,
$\B_{3}$ is $F^{\M}_{\l_{r}(\M )}$-constructible over $\B_{1}\cup\B_{2}$.

{\bf Claim 5.} For all $c\in\C_{1}$ there is
$B\subseteq\B_{1}$ of power $\le\l$ such that
$\tp (c,B)$ $F^{\M}_{\l^{+}}$-isolates $\tp (c,\B_{3})$.

\proof Assume not. Choose
$B_{1}\subseteq\B_{1}$ of power $\le\l$ and $c\in\C_{1}$
so that

(i) $\tp (c,B_{1})$ $F^{\M}_{\l^{+}}$-isolates
$\tp (c,\B_{1})$,

(ii) $\tp (c,B_{1})$ does not $F^{\M}_{\l^{+}}$-isolate
$\tp (c,\B_{3})$.

\noindent
By (ii) above, choose
$b\in\B_{3}$, $B_{0}\subseteq\B_{0}$
$C_{1}\subseteq\B_{1}$ and $C_{2}\subseteq\B_{2}$

(iii) $\vert C_{1}\cup C_{2}\vert <\l_{r}(\M )$,

(iv) $\tp (b,C_{1}\cup C_{2})$ $F^{\M}_{\l_{r}(\M )}$-isolates
$\tp (b,\B_{1}\cup\B_{2})$,

(v) $\tp (c,B_{1})$ does not $F^{\M}_{\l^{+}}$-isolate
$\tp (c,B_{1}\cup C_{1}\cup C_{2}\cup b)$,

(vi) for all $a\in B_{1}\cup C_{1}$, $\tp (a,\B_{2})$ does not split strongly
over $B_{0}$ and
$\vert B_{0}\vert\le\l$.

\noindent
Since $\B_{0}$ is $F^{\M}_{\l^{+}}$-saturated and $\M$ is $\l$-stable,
we can find $f\in Aut(B_{0})$
such that $f(C_{2})\subseteq\B_{0}$ and for all $a\in C_{2}$,
$f(a)\ E^{m}_{min,B^{*}_{0}}\ a$. Then by (vi),
$\tp (f(C_{2}),B_{0}\cup B_{1}\cup C_{1})=\tp (C_{2},B_{0}\cup B_{1}\cup C_{1})$.
Choose $b'\in\B_{1}$ so that
$\tp (b'\cup f(C_{2}),C_{1})=\tp (b\cup C_{2},C_{1})$. Then by (iv),
$\tp (b'\cup f(C_{2}),B_{0}\cup B_{1}\cup C_{1})=
\tp (b\cup C_{2},B_{0}\cup B_{1}\cup C_{1})$.
By (vi) $\tp (c,B_{1})$ does not $F^{\M}_{\l^{+}}$-isolate
$\tp (c,B_{1}\cup\C_{1}\cup f(C_{2})\cup b')$. Clearly this contradicts (i).
$\eop$ Claim 5

So $\B_{3}\cup\C_{1}$ is $F^{\M}_{\l_{r}(\M )}$-constructible over
$\C_{1}\cup\B_{2}$.

{\bf Claim 6.} For all $c\in\C_{2}$ there is
$B\subseteq\B_{2}$ of power $\le\l$ such that
$\tp (c,B)$ $F^{\M}_{\l^{+}}$-isolates $\tp (c,\C_{1}\cup\B_{3})$.

\proof As Claim 5 above. $\eop$ Claim 6.

So $\B_{3}\cup\C_{1}\cup C_{2}$ is $F^{\M}_{\l_{r}(\M )}$-constructible
(and so $F^{\M}_{\l^{+}}$-constructible) over
$\C_{1}\cup\C_{2}$. By the choice of $\C_{3}$, this implies that
$\C_{3}$ is
$F^{\M}_{\l^{+}}$-primary over
$\C_{1}\cup\C_{2}$.
$\eop$

Notice that in Theorem 6.3 the assumption, $\M$ is $\l$-stable,
is not necessary. We can avoid the use of it by Lemma 3.15.

\th 6.4 Lemma. Assume $\k >\l\ge\l (\M )$. Then $\l$-dop
implies $\k$-dop.

\proof Let $\A_{i}$, $i<4$, and $a$ as in the definition of
$\l$-dop. Choose $F^{\M}_{\k}$-saturated $\B_{0}\supseteq\A_{0}$
such that $\B_{0}\da_{\A_{0}}\A_{1}\cup\A_{2}$. Let $\B_{1}$ be
$F^{\M}_{\k}$-primary over $\B_{0}\cup\A_{1}$,
$\B_{2}$ be $F^{\M}_{\k}$-primary over $\B_{0}\cup\A_{2}$ and
$\B_{3}$ be $F^{\M}_{\k}$-primary over $\B_{1}\cup\B_{2}$.
Clearly we can choose the sets so that $\A_{3}\subseteq\B_{3}$
and $a\da_{\A_{3}}\B_{3}$.
By Lemmas 5.4 (iii) and 3.8 (iv),
$\B_{1}\da_{\B_{0}}\A_{2}$. Then
$\A_{2}\da_{\A_{0}}\B_{1}$ and so $\A_{2}\da_{\A_{1}}\B_{1}$.
By Lemma 5.4 (iii),

(1) $\A_{3}\da_{\A_{1}}\B_{1}$.

\noindent
Similarly,

(2) $\A_{3}\da_{\A_{2}}\B_{2}$.

\noindent
Also by Lemmas 5.4 (iii) and 3.8 (iv),
$\B_{1}\da_{\B_{0}}\B_{2}$.

By (1), (2), Lemma 4.5 and Corollary 4.8, $\tp (a,\B_{3})$ is orthogonal to
$\B_{1}$ and to $\B_{2}$. $\eop$

\th 6.5 Corollary. $\l (\M )$-dop implies $\l_{r}(\M )^{+}$-sdop.

\proof Immediate by Lemma 6.4 and Theorem 6.3. $\eop$

We finish this paper by giving open problems:

\th 6.6 Question. What are the relationships among the following
properties:

(1) $a\da_{A}A$,

(2) $a\nda_{A}a$,

(3) $t(a,A)$ is unbounded?

Notice that (1) does not imply (2) nor (3) (fails already
in the 'classical' case), (3) implies (2) (Lemma 3.2 (v))
and (1)$\wedge$(2) implies (3) (just choose $a_{i}$, $i<\vert\M\vert$,
so that $t(a_{i},A)=t(a,A)$ and
$a_{i}\da_{A}\cup_{j<i}a_{j}$).

\th 6.7 Question. Does Corollary 4.8 hold without the assumption
that the sets are strongly $F^{\M}_{\k_{r}(\M )}$-saturated?

\th 6.8 Question. Does the following hold: If $\M$ is superstable, then
for all $A$ there exists an '$a$-primary' set over $A$?

\vfill
\eject

\chapter{References.}

\item{[Hy]} T. Hyttinen, On nonstructure of elementary submodels
of a stable homogeneous structure, Fundamenta Mathematicae, vol. 156,
1998, 167-182.

\item{[Sh1]} S. Shelah, Finite diagrams stable in power, Annals
of Mathematical Logic, vol. 2, 1970, 69-118.

\item{[Sh2]} S. Shelah, Classification Theory, Stud. Logic Found. Math.
92, North-Holland, Amsterdam, 2nd rev. ed., 1990.

\bigskip

Tapani Hyttinen

Department of Mathematics

P.O. Box 4

00014 University of Helsinki

Finland

\medskip

Saharon Shelah

Institute of Mathematics

The Hebrew University

Jerusalem

Israel

\medskip

Rutgers University

Hill Ctr-Bush

New Brunswick

New Jersey 08903

U.S.A.

\end